\documentstyle[euscript,amsfonts,amssymb]{article}
\newcommand{\reals}{\Bbb R}
\newcommand{\rationals}{\Bbb Q}
\newcommand{\naturals}{\Bbb N}

\input mssymb

\newcommand{\QED}{\vrule width 6pt height 6pt depth 0pt \vspace{0.1in}}
\newcommand{\forces}{\mathrel{\|}\joinrel\mathrel{-}}

\newcommand{\JR}{J_{\reals}}
\newcommand{\ropo}{\rightarrow (\omega - path )^2_{\omega_1/<\omega }}
\newcommand{\Pk}{{\bold P}_{\kappa}}
\newcommand{\twoal}{\frak c}
\newcommand{\ropl}{\rightarrow (\omega - path )^2_{\lambda /<\omega }}
\newcommand{\VPk}{V^{\Pk}}
\newcommand{\pab}{p_{\alpha , \beta}}
\newcommand{\pag}{p_{\alpha , \gamma}}
\newcommand{\pbg}{p_{\beta , \gamma}}
\newcommand{\ab}{\{ \alpha , \beta\} }
\newcommand{\dKi}{{\dot K}_i}
\newcommand{\bpab}{{\bar p}_{\alpha , \beta }}
\newcommand{\bpz}{{\bar p}^0}
\newcommand{\bpo}{{\bar p}^1}
\newcommand{\bpaz}{\bpz_\alpha}
\newcommand{\bpgo}{\bpo_\gamma}
\newcommand{\bp}{\bar p}
\newcommand{\dX}{\dot X}
\newcommand{\dY}{\dot Y}
\newcommand{\dW}{\dot W}
\newcommand{\dZ}{\dot Z}
\def\JA{\langle J \rangle}
\def\JAR{\langle \JR \rangle }
\def\mg{MG(J)}
\def\smg{SMG(J)}
\def\vsg{VSG(J)}

\newtheorem{theorem}{Theorem}
\newtheorem{lemma}[theorem]{Lemma}
\newtheorem{prop}[theorem]{Proposition}
\newtheorem{corollary}[theorem]{Corollary}

\newtheorem{problem}{Problem}
\newtheorem{conjecture}{Conjecture}
\newtheorem{definition}{Definition}

\begin{document}

\title{Covering games and the Banach-Mazur game: $k$-tactics.}
\author{Tomek Bartoszynski\thanks{Supported by Idaho State Board of
   Education grant 92-096}\\ 
   Department of Mathematics,\\Boise State University,\\Boise, Idaho 83725
\and
    Winfried Just\thanks{Supported by NSF grant DMS-9016021
    and Research Challenge Grant RC 89-64 from Ohio University}\\Department of
    Mathematics,\\ Ohio University,\\Athens, Ohio 45701\\
\and 
    Marion Scheepers\thanks{Supported by Idaho State Board of Education grant
    91-093.}\\Department of Mathematics,\\Boise State University,\\Boise, Idaho
    83725} 
\date{}

\maketitle

\abstract{Given a free ideal $J$ of subsets of a set $X$, we consider  
   games where player ONE plays an increasing sequence of elements of the
   $\sigma$-completion of $J$, and player TWO tries to cover the union of this
   sequence by playing one set at a time from $J$. We describe various
   conditions under which player TWO has a winning strategy that uses
   only information about the most recent $k$ moves of ONE, and apply
   some of these results to the Banach-Mazur game.}

\bigskip

\section{Introduction}

Let $J$ be a free ideal of subsets of a given set. By $\JA$ we denote the
$\sigma$-ideal generated by $J$ ($\JA$ could turn out to be the power
set of $\cup J$). Two concrete examples of ideals motivated much of our work. 
The one is $\JR$, the ideal of nowhere dense subsets of the real line $\reals$.
In this case $\JAR$ is the ideal of meager sets of reals. The other is
$[\kappa]^{<\lambda}$ where $\omega=cof(\lambda)\leq\lambda\leq\kappa$
are cardinal numbers.

We are interested in games of the following type:
Player ONE plays a set $O_n \in \JA$ during inning $n$, to which TWO responds 
with a set $T_n \in J$. ONE is required to play an increasing sequence of
sets; TWO's objective is to cover $\bigcup_{n \in \omega} O_n$ with
$\bigcup_{n \in \omega } T_n$. As long as TWO remembers the complete history
of the game, this task is trivial. However, it often happens that
TWO needs to know only the last $k$ moves of the opponent in order to win.
A strategy that accomplishes this is called a {\it winning $k$-tactic}.

We consider four such games, $MG({\cal A},J)$, $\mg$, the ``monotonic game'',
$\smg$, the ``strongly monotonic game'', and $\vsg$, the ``very strong
game''. The study of these games was initiated in \cite{S1}, and motivated
by Telgarsk\'y's conjecture that for  
every $k > 0$ there exists a topological space $(X,\tau )$ such that
TWO has a winning $k+1$-tactic but no winning $k$-tactic in the Banach-Mazur
game on $(X,\tau )$ (see section 4.4 for more information).
However, we find the games considered here of interest independent 
of the original motivation. The game $\mg$ was introduced in
\cite{S1}, as was the game $\smg$; the games $MG({\cal A},J)$ and
$\vsg$ appear here for the first time. 

\smallskip

In sections 2 and 3, we introduce and discuss pseudo Lusin sets, the
irredundancy property and the coherent decomposition property of ideals. These
properties, together with the $\omega$-path partition relation, are
the main tools for constructing winning $k$-tactics in our games.
These combinatorial properties of ideals are very likely of independent
interest - they have already appeared in the literature in various guises.

In section 4 we apply the results of sections 2 and 3 to give various 
conditions sufficient for the existence of winning $k$-tactics for TWO in 
the games mentioned above. Not surprisingly, as the game becomes more 
favorable for TWO, weaker conditions suffice. 
Among other things, our results show that in the Banach-Mazur game on the 
space that inspired the invention of meager-nowhere dense games, TWO has
a winning 2-tactic.

The appendix is devoted to a proof of an unpublished consistency result
of Stevo Todorcevic, which we use in section 4.

\smallskip

   Our notation is mostly standard. One important exception may be that we
   use the symbol $\subset$ exclusively to mean ``is a proper subset of". Where
   we otherwise deviate from standard notation or terminology we explicitly
   alert the reader. For convenience we also assume the consistency of
   traditional (Zermelo-Fraenkel) set theory. All statements we make about the
   consistency of various mathematical assertions must be understood as
   consistency which can be proven by means of that theory. The reader
   might find having a copy of \cite{S1} and \cite{S2} handy when reading
   parts of this paper a bit more comfortable than otherwise.

   We are grateful to Stevo Todor\v{c}evi\'c for sharing with us his
   insights about the matters we study here, and for his kind permission to
   present in this paper some of his answers to our questions.

\section{The irredundancy property.}

      For a partially ordered set $(P,<)$ which has no maximum element we let
$$add(P,<)$$ be the least
   cardinal number, $\lambda$, for which there is a collection of cardinality
   $\lambda$ of elements of $P$ which do not have an upper bound in $P$. This
   cardinal number is said to be the {\em additivity} of $(P,<)$. Note
   that $add(P,<)$ is either $2$, or else it is infinite. In the latter case
   $(P,<)$ is said to be {\em directed}. We attend exclusively to directed
   partially ordered sets in this paper. Isbell \cite{I} and some earlier
   authors also refer to the additivity of a partially ordered set as its
   {\em lower character}; they denote it by $\ell(P,<)$.

   A free ideal $J$ on a set $S$ is partially ordered by $\subset$. The
   partially ordered set $(J,\subset)$ is directed. When
   $add(J,\subset)=\aleph_{0}$, the symbol $\langle J
   \rangle$ denotes the $\sigma$-completion of $J$ (i.e., the smallest
   collection which contains each union of countably many sets from $J$). We
   say that $J$ is a $\sigma$-{\em complete ideal} if $J=\langle J\rangle$.

   The other important example for our study is the set $^{\omega}\omega$ of
   sequences of nonnegative integers; we use ${\frak c}$ to denote the
   cardinality of this set. We say $g$ {\em eventually dominates}
   $f$ and write $f\ll g$ if: $\lim_{n\rightarrow\infty} (g(n)-f(n))=\infty$.
   It is customary to denote $add(^{\omega}\omega, \ll)$ by ${\frak b}$. 

   A well known theorem of Miller (\cite{M}, p. 94, Theorem 1.2) states that
   $$add(\langle J_{\reals} \rangle,\subset) \leq add(^{\omega}\omega,\ll)
   (={\frak b}).$$ 

   Again, for an arbitrary partially ordered set $(P,<)$ the symbol
   $$cof(P,<)$$ denotes the least cardinal number, $\kappa$, for which there is
   a collection $X$ of cardinality $\kappa$ of elements of $P$ such that: for
   each $p\in P$ there is an $x\in X$ such that $p\leq x$. This cardinal
   number is said to be the {\em cofinality} of $(P,<)$. Some authors (see
   e.g. \cite{I}, p. 397) also call this cardinal number the {\em upper
   character} of $(P,<)$ and denote it by $u(P,<)$. It is customary to
   denote $cof(^{\omega}\omega, \ll)$ by ${\frak d}$. 

   A theorem of Fremlin (\cite{F}, Proposition 13(b)) states that $$({\frak
   d}=) cof(^{\omega}\omega, \ll)\leq cof(\langle J_{\reals}
   \rangle,\subset).$$ 

   Let $(P,<)$ be a directed partially ordered set. The {\em bursting number}
   of $(P,<)$ (\cite{I}, p. 401) is the smallest cardinal number which
   exceeds the cardinality of each of the bounded subsets of $(P,<)$. This
   cardinal number is denoted by $burst(P,<)$. More important is the {\em
   principal bursting number} of $(P,<)$, denoted $bu(P,<)$ and
   define as $$bu(P,<)=\min\{burst(Q,<): Q\mbox{ is a cofinal subset of }P\}$$
   (following \cite{I}, p. 409). It is always the case that
   $add(P,<)\leq bu(P,<)$.

\begin{definition}
  A directed partially ordered set $(P,<)$ has the {\em irredundancy property}
  if: $$bu(P,<)= add(P,<).$$

  The cofinal subfamily ${\cal A}$ of $(P,<)$ is said to be {\em irredundant}
  if $burst({\cal A},<)\leq add(P,<)$. 
\end{definition}

   Not all $\sigma$-complete ideals have the irredundancy property. Here is
   an ad hoc example. Let $S_{1}$ and $S_{2}$ be disjoint sets such that
   $S_{i}$ has cardinality $\aleph_{i}$ for each $i$. Define an ideal $J$ on
   the union of these sets by admitting a set $Y$ into $J$ if: $Y\cap S_{1}$ is
   countable and $Y\cap S_{2}$ has cardinality less than $\aleph_{2}$. Then
   $add(J,\subset)=\aleph_{1}$ and $cof(J,\subset)=\aleph_{2}$. No cofinal
   family of $J$ is irredundant.

    A refined version of the classical notion of a Lusin set is
   instrumental in verifying the presence of the irredundancy property in many
   directed partially ordered sets. Since what we'll define is not exactly
   the same as the classical notion, we call our ``Lusin sets" {\em
   pseudo Lusin sets} (more about this after the definition). Let $\kappa$ and
   $\lambda$ be infinite cardinal numbers. Let $(P,<)$ be a directed partially
   ordered set.

\begin{definition} A subset $L$ of $P$ is a $(\kappa,\lambda)$
   pseudo Lusin set if:
\begin{enumerate}
\item{$\lambda$ is the cardinality of $L$ and}
\item{for each $x\in P$ the cardinality of the set $\{y\in L:y\leq x\}$
   is less than $\kappa$.} 
\end{enumerate} 
\end{definition}

  $(\kappa,\lambda)$ pseudo Lusin sets are interesting only when 
  $\kappa\leq\lambda$. If a directed partially ordered set
  $(P,<)$ has a 
  $(\kappa,\lambda)$ pseudo Lusin set, then $add(P,<)\leq\kappa$ and
  $\lambda\leq cof(P,<)$. Moreover, every partially ordered set
  has an $(add(P,<),add(P,<))$ pseudo Lusin set. Thus, if
  $add(P,<)=cof(P,<)$, then these are the only types of pseudo Lusin
  sets in $(P,<)$.

   Let $J$ be a free ideal on a set $S$. The uniformity number of $J$, written
   $unif(J)$, is the minimal cardinal $\kappa$ such that there is a
   subset of $S$ which is of cardinality $\kappa$, which is not an element
   of $J$.

   Consider the partially ordered set $(\langle
   J_{\reals}\rangle,\subset)$. If $L\subset {\reals}$ is a Lusin set in the
   classical sense (i.e., $L$ is uncountable and every meager set meets $L$ in
   only countably many points), then $\{\{x\}:x\in L\}$ is an $(\omega_1,\mid
   L\mid)$ pseudo Lusin set. There will be pseudo Lusin sets
   even when there are no (classical) Lusin sets: If
   $unif(\langle J_{\reals}\rangle)>add(\langle
   J_{\reals}\rangle,\subset)$ then every set of real numbers of cardinality
   $\aleph_1$ is meager, whence there is no Lusin set in the
   classical sense. Now let $\{M_{\alpha}:\alpha<add(\langle
   J_{\reals}\rangle,\subset)\}$ be a family of meager sets such that
\begin{enumerate}
\item{ $M_{\alpha}\subset M_{\beta}$ whenever $\alpha<\beta<add(\langle
   J_{\reals}\rangle,\subset)$ and}
\item{$\cup_{\alpha<add(\langle
   J_{\reals}\rangle,\subset)}M_{\alpha}$ is not meager.}
\end{enumerate}

   Then the set $L=\{M_{\alpha}:\alpha<add(\langle
   J_{\reals}\rangle,\subset)\}$ is a $(add(\langle
   J_{\reals}\rangle,\subset),add(\langle 
   J_{\reals}\rangle,\subset))$ pseudo Lusin set.
   
   It is also well known that these hypotheses on the ideal of meager
   subsets of the real line are consistent. For example, it is
   consistent that the real line is a union of $\aleph_1$ meager sets
   and that each set of real numbers of cardinality less than $\aleph_2$ is
   meager (see e.g. \cite{M}, \S6).

   The reader should also compare our notion of a $(\kappa,\lambda)$ -
   pseudo Lusin set with Cichon's notion of a $(\kappa,\lambda)$ - Lusin
   set (see \cite{Ci}).

   The connection between the irredundancy property and the existence
   of certain pseudo Lusin sets is given by the following proposition.
   The argument in its proof is well known in
   the special case when $P$ is the collection of
   countable subsets of some infinite set, ordered by set inclusion (see the
   proof of 4.4 on p. 409 of \cite{I}).

\begin{prop}\label{refinelusin}Let $(P,<)$ be a directed partially ordered
   set. Then the following statements are equivalent:
\begin{enumerate}
\item{There is an $(add(P,<),cof(P,<))$ pseudo Lusin set for $(P,<)$,}
\item{$(P,<)$ has the irredundancy property,}
\item{There is a cofinal $(add(P,<),cof(P,<))$ pseudo Lusin
   set for $(P,<)$,} 
\end{enumerate}
\end{prop}
\begin{description}\item[Proof.]{That 1. implies 2:\\
 Let \mbox{$L=\{x_{\xi}:\xi<cof(P,<)\}$} be
   such a pseudo Lusin set and let $\{a_{\xi}:\xi<cof(P,<)\}$ be a
   cofinal subfamily of $P$. For each $\xi<cof(P,<)$ choose $z_{\xi}\in P$
   such that $x_{\xi},a_{\xi}\leq z_{\xi}$. Put ${\cal
   A}=\{z_{\xi}:\xi<cof(P,<)\}$. Then ${\cal A}$ is an irredundant cofinal
   family.\\

That 2. implies 3:\\
   Let ${\cal A}$ be an irredundant cofinal family. We may assume that
the cardinality of this family is $cof(P,<)$. Then ${\cal A}$ is an
example of a cofinal $(add(P,<),cof(P,<))$ pseudo Lusin set.\\

It is clear that 3. implies 1.  $\QED$}
\end{description}

\begin{corollary}\label{cor:cardinalityideals} Let
$\kappa>\lambda\geq\aleph_0$ be cardinals, $\lambda$ regular. If
$cof([\kappa]^{<\lambda},\subset) = \kappa$, then
$([\kappa]^{<\lambda},\subset)$ has the irredundancy property.
\end{corollary}
\begin{description}\item[Proof.]{Let $\{S_{\alpha}:\alpha<\kappa\}$ be
a pairwise disjoint subcollection from $[\kappa]^{<\lambda}$. Then
this family is a $(\lambda,\kappa)$ pseudo Lusin set for this ideal.
Applying the cofinality hypothesis 
we conclude that this ideal has the irredundancy property. $\QED$}
\end{description}

   The ideal of finite subsets of an infinite set has the irredundancy
   property; the set of one-element subsets of such an infinite set forms an
   appropriate pseudo Lusin set for this ideal.

\begin{lemma}\label{equivalence} Let $\kappa>\lambda$ be an
uncountable cardinal numbers, $\lambda$ regular. Then the following
statements are equivalent: 
\begin{enumerate}
\item{ The ideal $([\kappa]^{<\lambda},\subset)$ has cofinality $\kappa$.}
\item{ There is a free ideal $J$ such that:
\begin{enumerate}
\item{$add(J,\subset)=\lambda$,}
\item{$cof(J,\subset)=\kappa$ and}
\item{$(J,\subset)$ has the irredundancy property.}
\end{enumerate}}
\end{enumerate}
\end{lemma}
\begin{description}\item[Proof.]{The proof of {\em 1}$\Rightarrow$ {\em 2} is
   trivial. We show that {\em 2} implies {\em 1}. Let $J$ be a
   free ideal on the set $S$ such that
   $cof(J,\subset)=\kappa$ and $add(J,\subset)=\lambda$, and
   $(J,\subset)$ has the irredundancy property. Let $L\subset J$ be an
   $(\lambda,\kappa)$ pseudo Lusin set for $J$. 
   Also let ${\cal C}\subset J$ be a cofinal family of cardinality $\kappa$.
   For each $X\in{\cal C}$ define: $S_X=\{Y\in L:Y\subseteq X\}$. Then the
   collection ${\cal B}=\{S_X: X\in {\cal C}\}$ is cofinal in
   $([L]^{<\lambda},\subset)$. $\QED$}
\end{description}

   The following examples play an important role in our
   game-theoretic applications.

\begin{center}{\bf {\large Example 1:} The ideal of countable subsets of an
   infinite set.}
\end{center}

   Let $\kappa$ be an uncountable cardinal number. Then
   $add([\kappa]^{\leq\aleph_{0}},\subset)=\aleph_{1}$ and 
   $bu([\kappa]^{\leq\aleph_{0}},\subset)\geq\aleph_{1}.$
   For uncountable cardinal numbers $\kappa$ it is always the case that
   $\kappa\leq cof([\kappa]^{\leq\aleph_{0}}, \subset)$.
   A set of the form $\{\{\alpha_{\xi}\}:\xi<\kappa\}$ (where this
   enumeration is bijective and $\lambda\leq\kappa$) is an
   $(\omega_{1},\kappa)$ pseudo Lusin set for $[\kappa]^{\leq\aleph_{0}}$. 
   The only difficult cases to decide whether or not the irredundancy
   property is present are those where 
   $\kappa< cof([\kappa]^{\leq\aleph_{0}}, \subset)$; this occurs for
   example when $\kappa$ has countable cofinality.
   It turns out that for these the irredundancy property is not 
   decidable by the axioms of traditional set theory:
\begin{enumerate}
\item{In \cite{To3}, Todor\v{c}evi\'{c} shows that if for each uncountable
   cardinal $\lambda$ of countable cofinality the assertions
\begin{enumerate}
\item{$cof([\lambda]^{\leq\aleph_{0}},\subset)=\lambda^{+}$ and}
\item{$\square_{\lambda}$}
\end{enumerate}
   are true, then for each uncountable cardinal number $\kappa$ there is a
   cofinal family ${\cal K}\subset [\kappa]^{\aleph_{0}}$ such that $|\{A\cap
   X:X\in{\cal K}\}|\leq\aleph_{0}$ for any countable subset $A$ of $\kappa$.
   Such a family ${\cal K}$ is an example of an
   $(add([\kappa]^{\leq\aleph_{0}},\subset),
   cof([\kappa]^{\leq\aleph_{0}},\subset))$ pseudo Lusin set for
   $([\kappa]^{\leq\aleph_{0}},\subset)$. These particular examples of
   pseudo Lusin sets are called {\em cofinal Kurepa families}.
   Thus it is true in the constructible universe,
   ${\bold L}$ that $([\kappa]^{<\aleph_0},\subset)$ has the
   irredundancy property for each infinite $\kappa$.}

\item{One might ask if any hypotheses beyond $ZFC$ are necessary to obtain the
   conclusion that $([\kappa]^{\leq\aleph_{0}},\subset)$ has the irredundancy
   property. Todor\v{c}evi\'{c} has shown in \cite{To2} that for an infinite
   cardinal number $\kappa$ the following statements are equivalent:
\begin{enumerate}
\item{$bu([\kappa]^{\leq\aleph_{0}},\subset)=\aleph_{1}$.}
\item{$([\kappa]^{\leq\aleph_{0}},\subset)$ has the irredundancy property.}
\end{enumerate}
   He also noted (p. 843 of \cite{To4}) that the version
   $$(\aleph_{\omega+1},\aleph_{\omega})\rightarrow(\omega_{1},\omega)$$ 
   of Chang's Conjecture implies that 
   $\aleph_{1}<bu([\aleph_{\omega}]^{\leq\aleph_{0}}, \subset)$ (and thus
   this ideal does not have the irredundancy property). Now \cite{L-M-S}
   established the consistency of the above version of Chang's Conjecture
   modulo the consistency of the existence of a fairly large cardinal.}
\item{This takes care of uncountable cardinals of countable
   cofinality. What is the situation for those of uncountable cofinality?
   It is clear that $([\kappa]^{\leq\aleph_0},\subset)$ has the
   irredundancy property if $\kappa$ is $\aleph_n$ for some finite $n$ or
   if, for some $m<\omega$, $\kappa$ is the $m$-th successor of a
   singular strong limit cardinal of uncountable cofinality. 
   In fact, the axiomatic system of traditional set theory has to be
   strengthened fairly dramatically before one could create circumstances
   where there is a 
   cardinal number of uncountable cofinality which is strictly less than the
   cofinality of its ideal of countable sets; it follows from Lemma 4.10 of
    \cite{J-M-P-S} that if there is a cardinal number of uncountable
    cofinality which is smaller than the cofinality of its ideal of countable
    sets, then there is an inner model with many measurable cardinal numbers. }
\end{enumerate}

    Information about the ideal of countable
   subsets of some infinite set can be used to gain information about
   some other ideals, using the notion of a locally small family. 

   \begin{definition}
  A family ${\cal F}$ of subsets of a set $S$ is {\em locally small}
  if:$$|\{Y\in {\cal F}:Y\subseteq X\}|\leq \aleph_{0}$$for each $X$ in ${\cal
  F}$. 
\end{definition}

  If the ideal of countable subsets of an infinite set has an
  irredundant cofinal family then that cofinal family is ipso facto locally
  small. If there is an
   $(\omega_{1},cof(J,\subset))$ pseudo Lusin set for the
   $\sigma$-complete free ideal $J$ on the set $S$, then $J$ contains
   a locally small cofinal family.

\begin{center}{\bf {\large Example 2:} The ideal of meager subsets of the
   real line}
\end{center}

   Assume that $add(\langle J_{\reals}\rangle,\subset)=cof(\langle
   J_{\reals}\rangle,\subset)$ (This equation is for example implied by Martin's
   Axiom). Then $\langle J_{\reals} \rangle$ has the irredundancy property. In
   this case one may insure that the cofinal family which witnesses the
   irredundancy is a well-ordered chain of meager sets. By the results cited
   from \cite{M} and \cite{F}, the hypothesis implies that ${\frak
b}={\frak d}$. It 
   is well known that the reverse implication is not
   provable.

   Irredundancy does not require having a well ordered cofinal chain of
   meager sets. For let an initial ordinal be given. According to a theorem
   of Kunen (\cite{K}, p. 906, Theorem 3.18) it is consistent that the
   cardinality of the real line is regular and larger than that initial
   ordinal, and at the same time there is an
   $(\omega_{1},{\frak c})$ pseudo Lusin set. It follows that $\langle
   J_{\reals}\rangle$ has a locally small cofinal 
   family of cardinality ${\frak c}$. In particular, $\langle J_{\reals
   }\rangle$ has the irredundancy property. If the continuum is larger than
   $\aleph_1$ it also follows that this ideal has no cofinal well-ordered
   chain.

   Stevo Todor\v{c}evi\'{c} has informed us that it is also consistent, modulo
   the consistency of a form of Chang's Conjecture that $\langle J_{{\reals
   }}\rangle$ does not	have the irredundancy property. Actually, something
   apparently weaker than that form of Chang's Conjecture is used: we present
   this result of Todor\v{c}evi\'{c}'s in Theorem \ref{consistency1},
   which he kindly permitted us to include in this paper. 

\begin{theorem}[Todor\v{c}evi\'{c}]\label{consistency1} If
``{\em ZFC+$MA_{\aleph_1}$+ there is no Kurepa family in
   $[\aleph_{\omega}]^{\aleph_0}$ of cardinality larger than
   $\aleph_{\omega}$}" 
   is a consistent theory, then so is the theory ``{\em ZFC \ + \
   $bu(\langle\JR\rangle , \subset )>add(\langle\JR\rangle, \subset)
   = \aleph_1$}". 
\end{theorem}
\begin{description}\item[Proof]{Let ${\bold P}$ be the set of finite functions
   with domain a subset of $\aleph_{\omega}$ and range a subset of $\omega$ (in
   other words, ${\bold P}$ is the standard set for adding
$\aleph_\omega$ Cohen 
   reals). For $p$ and $q$ in ${\bold P}$ we write $p<q$ if $q\subset p$. For
   $D$ a countable subset of $\aleph_\omega$ we write ${\bold P}(D)$
for the set 
   of elements of ${\bold P}$ whose domains are subsets of $D$.

   Suppose we have a sequence $\{N_{\xi}:\xi<\theta\}$
   ($\theta>\aleph_{\omega}$) of ${\bold P}$-names for meager sets of
reals. Let
   $D_{\xi}\in[\aleph_\omega]^{\aleph_0}$ be the support of $N_\xi$ i.e.,
   $N_\xi\in {\bold V}^{{\bold P}(D_{\xi})}$. By the hypothesis of the
theorem and
   by Theorem 1 of \cite{To3} there is an uncountable set $A\subset\theta$
   such that $D=\cup_{\xi\in A}D_{\xi}$ is countable. Thus, $N_{\xi}\in{\bold
   V}^{{\bold P}(D)}$ for each $\xi\in A$. Since ${\bold P}(D)$ is
essentially the 
   poset for adding {\em one} Cohen real and since $MA_{\aleph_1}$
holds, ${\bold 
   V}^{{\bold P}(D)}\models``\cup_{\xi\in A}N_{\xi} \mbox{ is
meager}"$ (because   
   ${\bold V}^{{\bold P}(D)}\models``MA(\sigma-\mbox{centered})"$). $\QED$}
\end{description}

   The hypothesis of Theorem \ref{consistency1} is consistent modulo
   the consistency of the relevant form of Chang's Conjecture, because
   that form of the conjecture is preserved by c.c.c. generic extensions. 

\section{The coherent decomposition property}

Let $J$ be a free ideal on a set $S$ and let $\langle J\rangle$
be its $\sigma$-completion. Let ${\cal A}$ be a subcollection of $\langle
J\rangle$.

\begin{definition} 
\begin{enumerate}
\item{${\cal A}$ has a coherent decomposition if there is for each
$A\in {\cal A}$ a sequence $(A^n:n<\omega)$ such that:
\begin{enumerate}
\item{$A^n\in J$ for each $n$,}
\item{$A^n\subseteq A^m$ whenever $n<m<\omega$, and}
\item{For all $A$ and $B$ in ${\cal A}$ such that $A\subset B$, there is an
$m$ such that $A^n\subseteq B^n$ whenever $n\geq m$.}
\end{enumerate}
  The collection $\{(A^n:n<\omega): A\in {\cal A}\}$ is said to be a coherent
  decomposition for ${\cal A}$.}
\item{The ideal $J$ has the coherent decomposition property if some cofinal
subset of $\langle J\rangle$ has a coherent decomposition.}
\end{enumerate}
\end{definition}

It is worth mentioning that if $J$ has the coherent decomposition
property and if $\langle J\rangle$ has a cofinal chain, than the
family $\langle J\rangle$ itself has a coherent decomposition. 
We now explore the coherent decomposition property for our examples.

\begin{center}{\bf {\large Example 1:} (continued)}
\end{center}

\begin{theorem}\label{locsmallth} Let ${\cal A}$ be a locally small family of
    countable sets such that $({\cal A},\subset)$ is a well-founded partially
   ordered set. Then ${\cal A}$ has a coherent decomposition.
\end{theorem}

\begin{description}\item[Proof.]{ Let $\Phi:{\cal A}\rightarrow\alpha$ be a
   function to an ordinal $\alpha$ such that $\Phi(A)<\Phi(B)$ for all $A\subset
   B$ in ${\cal A}$ (i.e., a rank function). Since ${\cal A}$ is locally small we
   may assume that $\alpha$ is $\omega_{1}$.

   For $A$ in ${\cal A}$ with $\Phi(A)=0$, choose a sequence $(A^{n}:n<\omega)$
   of finite subsets of $A$ such that $A=\cup_{n<\omega}A^{n}$ and
   $A^{n}\subseteq A^{n+1}$ for all n.

   Let $0<\beta<\omega_{1}$ be given and assume that we have already assigned to
   each $A$ in ${\cal A}$ for which $\Phi(A)<\beta$, a sequence
   $(A^{n}:n<\omega)$ in compliance with {\em 1} and {\em 2}. Now Let $B$ be an
   element of ${\cal A}$ such that $\Phi(B)=\beta$. Write $F(B)=\{A\in {\cal
   A}: A\subset B\}.$
 
   To begin, arbitrarily choose a sequence $(S_{n}:n<\omega)$ of finite sets
   such that $B=\cup_{n<\omega}S_{n}$. For each
   $A\in F(B)$, define $g_{A}:\omega\rightarrow\omega$ such that for each
   $n<\omega$, $$g_{A}(n)=min\{k<\omega:A^{n}\subseteq S_{0}\cup\dots\cup
   S_{k}\}.$$ Then $\{g_{A}:A\in F(B)\}$ is countable since ${\cal A}$ is
   locally small. Let $f\in$ \mbox{$^{\omega}\omega$} be a strictly increasing
   function such that $g_{A}\ll f$ for each $A$ in $F(B)$. Define:
   $$B^{n}=S_{0}\cup\dots\cup S_{f(n)}$$ for each n. Then $(B^{n}:n<\omega)$ is
   as required. $\QED$}
\end{description}

\begin{corollary}\label{localsmallcor}Let $J$ be a free ideal on a set $S$ and
   let ${\cal A}$ be a locally small family of sets in $\langle J\rangle$ such
   that $({\cal A},\subset)$ is a well-founded partially ordered set. Then
   ${\cal A}$ has a coherent decomposition. 
\end{corollary}

\begin{description}\item[Proof.]{ For each $B$ in ${\cal A}$, let
   $(S_{n}(B):n<\omega)$ be a sequence from $J$ such that
   $B=\cup_{n<\omega}S_{n}(B)$. Also write $\Gamma(B)=\{A\in{\cal A}:A\subseteq
   B\}$. Then \mbox{${\cal B}=\{\Gamma(A):A\in {\cal A}\}$} is a well-founded,
   locally small collection of countable subsets of ${\cal A}$. Choose, by
   Theorem \ref{locsmallth}, for each $A\in {\cal A}$ a sequence
   $(\Gamma(A)^{n}:n<\omega)$ of finite subsets of $\Gamma(A)$ such that:
\begin{enumerate}
\item{$\Gamma(A)=\cup_{n<\omega}\Gamma(A)^{n}$ where $\Gamma(A)^{n}\subseteq
   \Gamma(A)^{n+1}$ for each n, and}
\item{for all $A$ and $B$ in ${\cal A}$ with $A\subset B$ there exists an $m$
   such that: $$\Gamma(A)^{n}\subseteq \Gamma(B)^{n}$$ for
   all $n\geq m$.}
\end{enumerate}

			For each $A$ in ${\cal A}$ and each $n<\omega$ define:
   $$A^{n}=\cup\{S_{j}(B):j\leq n \mbox{ and } B\in\Gamma(A)^{n}\}.$$ Then the
   sequences $(A^{n}:n<\omega)$ are as required. $\QED$}
\end{description}

\begin{corollary} If $([\kappa]^{\leq \aleph_0},\subset)$ has the
irredundancy property, then it has the coherent decomposition property.
\end{corollary} 

\begin{description}\item[Proof]{An irredundant cofinal family is necessarily
   locally small. We may thin out any cofinal family to a well-founded cofinal
   family. Now apply Theorem \ref{locsmallth}. $\QED$}
\end{description}  

\begin{center}{\bf {\large Example 2:} (continued)}
\end{center}

   We show that the ideal of meager sets of the real line
   has the coherent decomposition property, and also that it has a second
   combinatorial property which plays an important role in our
   game-theoretic applications. It is convenient, for this section, to
   work with the set $^{\omega}2$, with the usual Tychonoff product
   topology ($2=\{0,1\}$ is taken to have the discrete topology) in place
   of $\reals$. For a subset $S$ of the domain of a function $g$, the
   symbol $g\lceil_{S}$ denotes the restriction of $g$ to the set $S$.
   For $s$ an element of $^{<\omega}2$, the symbol $[s]$ 
   denotes the set of all those $x$ in $^{\omega}2$ for which
   $x\lceil_{length(s)}=s$. Subsets of $^{\omega}2$ of the form $[s]$
   where $s$ ranges over $^{<\omega}2$, form a base for the topology of
   $^{\omega}2$. 
   Let $f\in\mbox{$^{\omega}\omega$}$ be a strictly increasing sequence and let
   $x$ be an element of $^{\omega}2$. Define:

$$B_{x,f}=\{z\in\mbox{$^{\omega}2$}:\forall^{\infty}_n(z\lceil_{[f(n),f(n+1))}\not=
   x\lceil_{[f(n),f(n+1))})\}.
$$

   Now also fix an $n\in \omega$ and define

$$B^n_{x,f}= \{z\in \mbox{$^{\omega}2$}:(\forall k\geq
   n)(z\lceil_{[f(k),f(k+1))}\not=x\lceil_{[f(k),f(k+1))})\}.
$$
   Then $B^m_{x,f}\subseteq B^n_{x,f}$ whenever $m<n<\omega$; also,
   $B_{x,f}=\cup_{n<\omega}B^n_{x,f}$.

\begin{prop}\label{crucial} For $x,y\in\mbox{$^{\omega}2$}$ and strictly
   increasing $f,g\in \mbox{$^{\omega}\omega$}$, the following assertions are
   equivalent: 
\begin{enumerate}
\item{$B_{x,f}\subset B_{y,g}$.}
\item{\begin{enumerate}
   \item{$B_{x,f}\not=B_{y,g}$ and}
   \item{$(\forall^{\infty}_n) (\exists k) (g(n)\leq f(k)<f(k+1)\leq
         g(n+1)\mbox{ and }
         x\lceil_{[f(k),f(k+1))}=y\lceil_{[f(k),f(k+1))})$} 
   \end{enumerate}}
\end{enumerate}
\end{prop}

\begin{description}\item[Proof.]{

   That {\em 1} implies {\em 2} requires some thought:\\

   If {\em 1} holds, then (a) of {\em 2} holds. Assume
   the negation of {\em 2}(b). It reads:
$$(\exists^{\infty}_n)(\forall k)(\neg(g(n)\leq f(k)<f(k+1)\leq
   g(n+1)) \mbox{ or
   }\neg(x\lceil_{[f(k),f(k+1))}=y\lceil_{[f(k),f(k+1))}))
$$

   Put $S=\{n<\omega:(\forall
k)(\neg([f(k),f(k+1)]\subseteq[g(n),g(n+1)]) \mbox{ or }\neg(
x\lceil_{[f(k),f(k+1))}= y\lceil_{[f(k),f(k+1))})\}$. Our hypothesis
is that $S$ is an infinite set.

   Consider an $n$ in $S$. For each $k$, there are the
following possibilities:
\begin{enumerate}
\item{$\neg([f(k),f(k+1)]\subseteq[g(n),g(n+1)]$}
\item{$[f(k),f(k+1)]\subseteq[g(n),g(n+1)]$, but
$x\lceil_{[f(k),f(k+1))}\not= y\lceil_{[f(k),f(k+1))}$.}
\end{enumerate}

   Put $S_n=\{k:\mbox{{\em 2} holds for }k\}$.
We consider two cases.\\ 
  
   {\bf Case 1:} There are infinitely many $n$ for which $S_n$ is
nonempty.\\
   Choose an infinite sequence $(n_1,n_2,n_3,\dots)$ from $S$ such that:

\begin{enumerate}
\item{$S_{n_m}\neq\emptyset$,}
\item{$n_{m+1}>g(n_m+1)$, and}
\item{$(\exists k)(g(n_m+1)<f(k)<g(n_{m+1}))$, for each $m$, and}
\item{$f(1)<g(n_1)$.}
\end{enumerate}

   This is possible because $f$ and $g$
are increasing, and $S$ is infinite. Put
$T=\cup_{j=1}^{\infty}[g(n_j),g(n_j+1))$. Define $z$, an element
of $^{\omega}2$, so that $z\lceil_{T}=y\lceil_{T}$ and $z(n)=1-x(n)$
for each $n\in\omega\backslash T$. Then $z\in B_{x,f}$ while $z\not\in
B_{y,g}$. Thus {\em 1} fails in this case.\\

   {\bf Case 2}: There are only finitely many $n\in S$ for which
$S_n$ is nonempty.\\

   We may assume that $S_n=\emptyset$ for each $n\in S$. Consider
   $n\in S$. We then have that for each $k\in \omega$,
   $[f(k),f(k+1))\not\subseteq[g(n),g(n+1))$. We distinguish between
   two possibilities:
\begin{enumerate}
\item{$(\exists k)(g(n)\leq f(k)<g(n+1))$ or}
\item{$(\forall k)(f(k)\not\in[g(n),g(n+1))$}
\end{enumerate}

{\bf Case 2 (A):} Possibility {\em 1} occurs for infinitely many $n\in S$:\\

   Choose $n_1<n_2<n_3<\dots$ from $S$ such that
\begin{itemize}
\item{$2\cdot n_j\leq n_{j+1}$ for each $j$,}
\item{for each $j$ there is a $k$ such that
   $g(n_j+1)<f(k)<g(n_{j+1})$,}
\item{for each $j$ there is a $k$ such that
   $f(k)\in[g(n_j),g(n_{j+1}))$, and}
\item{$f(1)<g(n_1)$.}
\end{itemize}
   Put $T=\cup_{j=1}^{\infty}[g(n_j),g(n_j+1))$ and define $z$ so that
$z\lceil_T=y\lceil_T$, and $z(n)=1-x(n)$ for each
$n\in\omega\backslash T$. From the hypothesis of Case 2(A) it follows
that $z\in B_{x,f}$, but $z\not\in B_{y,g}$. Thus, {\em 1} of the
Proposition fails also in this case.\\

{\bf Case 2 (B):} Possibility {\em 1} occurs for only finitely many $n\in S$:\\

   We may assume that possibility {\em 2} occurs for each $n\in
   S$. Choose $k_1<k_2<k_3<\dots$ such that for each $j$ there is
   an $n\in S$ with $[g(n),g(n+1))\subset[f(k_j),f(k_j+1))$. For
   each $j$ choose $n_j\in S$ such that $[g(n_j),g(n_j+1))\subset
   [f(k_j),f(k_j+1))$. As before define
   $T=\cup_{j=1}^{\infty}[g(n_j),g(n_j+1))$. Finally, define $z$ so
   that $z\lceil_T=y\lceil_T$ and $z(n)=1-x(n)$ for each
   $n\in\omega\backslash T$. Then $z\in B_{x,f}$ and $z\not\in B_{y,g}$,
   showing that {\em 1} of the Proposition fails also in this case.

   This completes the proof of the Proposition. $\QED$}
\end{description}

\begin{lemma}\label{equality}Let $f$ and $g$ be strictly increasing
elements of $^{\omega}\omega$ for which there is some $k<\omega$ such
that $g(n+k)=f(n)$ for all but finitely many $n$. If $B_{x,f}\subseteq
B_{y,g}$, then $B_{x,f}=B_{y,g}$.
\end{lemma}

\begin{description}\item[Proof.]{ Assume that $B_{x,f}\neq B_{y,g}$
and suppose that $B_{y,g}\not\subseteq B_{x,f}$.
   We show that $B_{x,f}\not\subseteq B_{y,g}$. Let $z$ be an element of
   $B_{y,g}\backslash B_{x,f}$. Fix $N$ such that
\begin{enumerate}
\item{$z\lceil_{[g(n+k),g(n+k+1))}\neq y\lceil_{[g(n+k),g(n+k+1))}$
   and}
\item{$f(n)=g(n+k)$}
\end{enumerate}
   for each $n\geq N$.

   Since $z$ is not an element of $B_{x,f}$, there are infinitely many
   $n\geq N$ for which $z\lceil_{[f(n),f(n+1))}=x\lceil_{f(n),f(n+1))}$.
   Consequently the set  $S=\{n\geq N:x\lceil_{[f(n),f(n+1))}\neq
   y\lceil_{[f(n),f(n+1))}\}$ is infinite. Now define $t$
   such that $t\lceil_{[f(n),f(n+1))}=y\lceil_{[f(n),f(n+1))}$ for each
   $n\in S$, and $t(m)=1-x(m)$ for each $m\in\omega\backslash(\cup_{n\in
   S}[f(n),f(n+1)))$. Then $t$ is in $B_{x,f}$ but not in $B_{y,g}$.
   $\QED$}
\end{description}

   Under the hypothesis of Lemma \ref{equality},
   $x(n)=y(n)$ for all but finitely many $n$.

\begin{prop}\label{orderpres} Let $x,y$ be elements of $^{\omega}2$
   and let $f,g$ be increasing elements of $^{\omega}\omega$. Of the
   following two assertions, {\em 1} implies {\em 2}.
\begin{enumerate}
\item{$B_{x,f}\subset B_{y,g}$.}
\item{$f\ll g$.}
\end{enumerate}
\end{prop}

\begin{description}\item[Proof.]{ Assume that $B_{x,f}\subset
   B_{y,g}$. Fix, by Proposition \ref{crucial}, an $N$ such that\\

$(\forall n\geq N) (\exists k) ([f(k),f(k+1)]\subseteq[g(n),g(n+1)]$
   and $x\lceil_{[f(k),f(k+1))}=y\lceil_{[f(k),f(k+1))}).$\\

   For each $n\geq N$ choose $k_n$ such that
   $[f(k_n),f(k_n+1)]\subseteq [g(n),g(n+1)]$. It follows that $k_n+1\leq
   k_{n+1}$ for each $n\geq N$ (since $f$ and $g$ are increasing).\\
{\bf Claim:} $[f(k_n),f(k_n+1)]\subset [g(n),g(n+1)]$ for infinitely many n.\\
{\bf Proof of the claim:} For otherwise, fix $M\geq N$ such that
   $[f(k_n),f(k_{n+1}]=[g(n),g(n+1)]$ for each $n\geq M$. Then we have
   $k_{n+1}=k_n+1$ for each $n\geq M$. It follows that
   $g(n)=f(n+(k_M-M))$ for all $n\geq M$. Then Lemma \ref{equality}
   implies that $B_{x,f}=B_{y,g}$, contrary to the fact that $B_{x,f}$ is
   a proper subset of $B_{y,g}$. This completes the proof of the claim.
$\QED$\\

   Thus, there are infinitely many $n$ for which $k_{n+1}>k_n+1$. Let
   $m>1$ be given, and fix $L\geq M$ such that $|\{n<L:k_{n+1}>k_n+1\}|\geq
   k_1+m$. Then $k_n>(n+m)$ for each $n\geq L$; we have
$$f(n+1)<f(n+m)\leq f(k_n)<g(n+1)$$
   for each $n\geq L$. In particular, $m\leq g(n+1)-f(n+1)$ for each
   $n\geq L$. This completes the proof that $f\ll g$. $\QED$}
\end{description}

\begin{prop}\label{coherencyprop} Let $x$ and $y$ be elements of
   $^{\omega}2$ and let $f$ and $g$ be increasing elements of
   $^{\omega}\omega$. If $B_{x,f}\subset B_{y,g}$, then there is an
   $m<\omega$ such that $B^n_{x,f}\subseteq B^n_{y,g}$ whenever $n\geq m$.
\end{prop}

\begin{description}\item[Proof]{From our hypotheses and Proposition
\ref{crucial} there is an $m$ such that for each $n\geq m$
there is a $k$ such that $[f(k),f(k+1))\subseteq[g(n),g(n+1))$ and
$x\lceil_{[f(k),f(k+1))}=y\lceil_{[f(k),f(k+1))}$. By Proposition 
\ref{orderpres} there is an $M>m$ such that $f(j)\leq
g(j)$ for each $j\geq M$. We show that $B^n_{x,f}\subseteq B^n_{y,g}$
for each $n\geq M$. 

Let $z$ be an element of $B^n_{x,f}$. Then
$z\lceil_{[f(j),f(j+1))}\neq x\lceil_{[f(j),f(j+1))}$ for each $j\geq
n$. But consider any $j\geq $n. Then there is a $k$ such that
$[f(k),f(k+1))\subset [g(j),g(j+1))$; $k\geq j$ for any such $k$, by
the choice of $M$. It follows that $z\lceil_{[g(j),g(j+1))}\neq
y\lceil_{[g(j),g(j+1))}$. Thus, $z$ is also an element of $B^n_{y,g}$. $\QED$}
\end{description}

\begin{prop}\label{cofinalprop} For each $X\in \langle
   J_{\reals}\rangle$ there are an $x$ in $^{\omega}2$ and an increasing
   $f$ in $^{\omega}\omega$ such that $X\subset B_{x,f}$.
\end{prop}
\begin{description}\item[Proof.]{ Let $X$ be a meager set. We may
   assume that $X=\cup_{n=0}^{\infty}X_n$ where $X_n\subseteq X_{n+1}$ and
   $X_n$ is closed, nowhere dense for each $n$. Fix a well-ordering
   of $^{<\omega}2$, and define $(s_n:n<\omega)$ and $f$ in
   $^{\omega}\omega$ as follows:\\
   Take $s_0=\emptyset$ and $f(0)=0$. Assume that
   $s_1,s_2,\dots,s_n$ and $f(1),\dots,f(n)$ have been defined so that:
\begin{enumerate}
\item{$s_1$ is the first element of $^{<\omega}2$ such that $[s_1]\cap
   X_1=\emptyset$ and $f(1)=length(s_1)$,}
\item{$s_{j+1}$ is the first element of $^{<\omega}2$ such that
   $[t^{\frown}s_{j+1}]\cap X_{j}=\emptyset$ for each $t$ in $^{\leq
   f(j)}2$, and $f(j+1)=\sum_{i=0}^{j+1}length(s_i)$ for each $j<n$.}
\end{enumerate}

   Then let $s_{n+1}$ be the first element of $^{<\omega}2$ such that
   $[t^{\frown}s_{n+1}]\cap X_n=\emptyset$ for each $t$ in $^{\leq
   f(n)}2$; put $f(n+1)=f(n)+length(s_{n+1})$. 

   Finally, set $x=s_1^{\frown}s_2^{\frown}s_3^{\frown}\dots$.\\
{\bf Claim:} $X\subseteq B_{x,f}$.\\
   For suppose that $z$ is not an element of $B_{x,f}$. Then there are
   infinitely many $n$ for which
   $z\lceil_{[f(n),f(n+1))}=x\lceil_{[f(n),f(n+1))}$; in other words,
   there are infinitely many $n$ for which
   $z\lceil_{[f(n),f(n+1))}=s_{n+1}$. Now fix an $m$. Choose an
   $n>m$ such that $z\lceil_{[f(n),f(n+1))}=s_{n+1}$. From the choice of
   $s_{n+1}$ it follows that $[z\lceil_{f(n+1)}]\cap X_m=\emptyset$; in
   particular, $z\not\in X_m$. Consequently, $z$ is not an element
   of $X$. $\QED$}
\end{description}

\begin{prop}\label{meagernessprop} Each $B^n_{x,f}$ is in $J_{\reals}$.
\end{prop}
\begin{description}\item[Proof.]{ Consider an $s$ from $^{\omega}2$
   for which $[s]\cap B^n_{x,f}\neq\emptyset$. Choose $m$ such that
   $f(m)>length(s)$ and $m>n$. Then choose $t$ from $^{<\omega}2$ such that
   $length(s^{\frown}t)\geq f(m+1)$ and
   $s^{\frown}t\lceil_{[f(m),f(m+1))}= x\lceil_{[f(m),f(m+1))}$. Then
   $[s^{\frown}t]\cap B^n_{x,f}=\emptyset$. It follows that $B^n_{x,f}$
   is nowhere dense. $\QED$}
\end{description} 

   Consequently, $B_{x,f}$ is a meager set for each $x$ in
   $^{\omega}2$ and for each increasing $f$ from $^{\omega}\omega$.

\begin{theorem}\label{important} $\langle J_{\reals}\rangle$ has a
   cofinal family which embeds in $(^{\omega}\omega,\ll)$ and which has
   the coherent decomposition property.
\end{theorem}

\begin{description}\item[Proof.]{ By Propositions \ref{meagernessprop}
   and \ref{cofinalprop} the family of sets of the form $B_{x,f}$ where
   $f$ is an increasing element of $^{\omega}\omega$ and $x$ is an
   element of $^{\omega}2$, is a cofinal family of meager sets. By
   Proposition \ref{coherencyprop}, this family has the coherent
   decomposition property. Also, the mapping which assigns $f$ to
   $B_{x,f}$ is, according to Proposition \ref{orderpres}, an order
   preserving mapping. $\QED$}
\end{description}

\begin{center}{\bf {\large Example 3:} Cardinals of countable cofinality}
\end{center}

 Here is a result which is quite analogous to Theorem \ref{locsmallth}.

\begin{theorem}\label{unctble,countblecof} Let $\lambda$ be an uncountable
   cardinal number which has countable cofinality. Let
   $\lambda_{0}<\lambda_{1}<\dots$ be a sequence of infinite regular cardinal
   numbers which converges to $\lambda$. Let $({\cal A},\subset)$ be a
   well-founded family of sets, each of cardinality $\lambda$, such that 
   $$|\{Y\in{\cal A}:Y\subseteq X\}|\leq\lambda$$ for each $X$ in ${\cal A}$.
   Then ${\cal A}$ has the coherent decomposition property. In
   particular:\\
   There exists for each $A\in{\cal A}$ a sequence $(A^{n}:n<\omega)$ such
   that:
\begin{enumerate}
\item{$|A^{n}|\leq\lambda_{n}$ for all $n$,}
\item{$A^{n}\subseteq A^{n+1}$ for all n,}
\item{$A=\cup_{n=0}^{\infty}A^{n}$ and}
\item{if $A\subset B$, then there is an $m<\omega$ such that $A^{n}\subseteq
   B^{n}$ for all $n\geq m$.}
\end{enumerate}
\end{theorem}
\begin{description}\item[Proof.]{Let $\Phi:{\cal A}\rightarrow\lambda^{+}$ be a
   rank function. For all $A$ in ${\cal A}$ with $\Phi(A)=0$, choose
   $(A^{n}:n<\omega)$ arbitrary, subject only to {\em 1, 2} and {\em 3}.

   Let $0<\gamma<\lambda^{+}$ be given and assume that $(A^{n}:n<\omega)$ has
   been assigned to each $A$ from ${\cal A}$ for which $\Phi(A)<\gamma$, in such
   a way that {\em 1, 2, 3} and {\em 4} are satisfied. Consider $B$ in
   ${\cal A}$ with $\Phi(B)=\gamma$. Write $F(B)$ for $\{A\in{\cal A}:A\subseteq
   B\}$ and write $F(B)=\cup_{n=0}^{\infty}F_{n}(B)$ where 
\begin{enumerate}
\item{$F_{0}(B)\subseteq F_{1}(B)\subseteq\dots$, and}
\item{$|F_{n}(B)|\leq\lambda_{n}$ for all $n$.}
\end{enumerate}
  Also let $B=\cup_{n=0}^{\infty}X_{n}$ where $X_{0}\subseteq
  X_{1}\subseteq\dots$ and $X_{n}\leq \lambda_{n}$ for all $n$. Finally put
  $B^{n}=(\cup\{A^{n}:A\in F_n(B)\})\cup X_{n}$ for each $n$. Then
  $(B^{n}:n<\omega)$ is as required. $\QED$}
\end{description}

\begin{corollary}\label{cor:countablecof} Let $\lambda$ be a cardinal number
   of countable cofinality. If $([\kappa]^{\leq\lambda},\subset)$
   has the irredundancy property then it has the coherent
   decomposition property. 
\end{corollary}

\section{Applications}

   The $\omega-path$ partition relation is the one other combinatorial
   ingredient in our technique for constructing winning $k$-tactics, or
   for defeating a given $k$-tactic for TWO. For a
   positive integer $n$, infinite cardinal number $\lambda$ and a
   partially ordered set $(P,<)$, the symbol
$$
   (P,<)\rightarrow(\omega-path)^{n}_{\lambda/<\omega}
$$
   means that for every function $F:[P]^{n}\rightarrow\lambda$ there is an
   increasing $\omega$-sequence $$p_{1}<p_{2}<\dots<p_{m}<\dots$$ such that the
   set $\{F(\{p_{j+1},\dots,p_{j+n}\}):j<\omega\}$ is finite. 
   The negation of this assertion is denoted by the symbol 
   $$(P,<)\not\rightarrow(\omega-path)^{n}_{\lambda/<\omega}.$$
   This partition relation has been studied in \cite{S2}. The reader
   should consult this reference about the various facts concerning the
   $\omega$-path relation which are used in the sequel. 

\subsection{The game $MG({\cal A},J)$}
   
   For a free ideal $J$ on an infinite set $S$ and for a family ${\cal
   A}$ in $\langle J\rangle$ with the property that for 
   each $X\in{\cal A}$ there is a $Y\in{\cal A}$ such that $X\subset Y$, the
   game $MG({\cal A},J)$ is defined so that an $\omega$-sequence
   $(O_{1},T_{1},\dots,O_{n},T_{n},\dots)$ is a play if for each $n$,
\begin{enumerate}
\item{$O_{n}\in {\cal A}$ is player ONE's move in inning $n$,}
\item{$T_{n}\in J$ is player TWO's move in inning $n$, and}
\item{$O_{n}\subset O_{n+1}$.}
\end{enumerate}
   Player TWO wins this play if
   $\cup_{n=1}^{\infty}O_{n}\subseteq\cup_{n=1}^{\infty}T_{n}.$

\begin{theorem}\label{mgth} Let $J$ be a free ideal on a set $S$. If ${\cal A}$
   is a family of sets in $\langle J\rangle$ such that: 
\begin{enumerate}
\item{for each $X\in{\cal A}$ there is a $Y\in{\cal A}$ such that $X\subset Y$,}
\item{$({\cal A},\subset)\not\rightarrow(\omega-path)^{k}_{\omega/<\omega}$
   for some $k\geq 2$, and}  
\item{${\cal A}$ has a coherent decomposition}
\end{enumerate}
   then TWO has a winning $k$-tactic in $MG({\cal A},J)$.
\end{theorem}
\begin{description}\item[Proof.]{ Choose a function $F:[{\cal
   A}]^{k}\rightarrow\omega$ which witnesses hypothesis {\em 2}. Also
   associate with each $A$ in ${\cal A}$ a sequence $(A^{n}:n<\omega)$ such
   that hypothesis {\em 3} is satisfied. 
   
   Define a $k$-tactic, $\Upsilon$ for TWO as follows. Let
   $(X_{1},\dots,X_{j})$ be given such that $j\leq k$, 
   $X_{1}\subset\dots\subset X_{j}$ and $X_{i}\in{\cal A}$ for $i\leq j$.
\begin{enumerate}
\item{If $j<k$: Then put $\Upsilon(X_{1},\dots,X_{j}) =
   X^{1}_{1}\cup\dots\cup X^{1}_{j}.$} 
\item{If $j=k$: Let $m$ be such that
\begin{itemize}
\item{$m\geq F(\{X_{1},\dots,X_{k}\})$ and}
\item{$X^{n}_{1}\subseteq\dots\subseteq X^{n}_{k}$ for all $n\geq m$.}
\end{itemize}
    Put $\Upsilon(X_{1},\dots,X_{k})=X^{m}_{1}\cup\dots\cup X^{m}_{k}$.}
\end{enumerate}

  Then $\Upsilon$ is a winning $k$-tactic for TWO. For let 
  $(O_{1},T_{1},\dots,O_{n},T_{n},\dots)$ be a play of $MG({\cal A},J)$
  where:
\begin{itemize}
\item{$T_{j}=\Upsilon(O_{1},\dots,O_{j})$ for each $j\leq k$}
\item{$T_{n+k}=\Upsilon(O_{n+1},\dots,O_{n+k})$ for each $n<\omega$.}
\end{itemize}
   For each $t\geq 1$ let $m_{t}$ be the number associated with
   $(O_{t},\dots,O_{t+k-1})$ in part 2 of the definition of $\Upsilon.$
   By the properties of $F$, the set $\{m_{t}:t=1,2,3,\dots\}$ is
   infinite.
   Thus choose $t_{1}<t_{2}<\dots$ such that $m_{j}<m_{t_{r}}$ for all
   $j<t_{r}$. It follows from the criteria used in the choices of the
    numbers $m_{t}$ 
   that $$O^{m_{t_{r}}}_{1}\subseteq\dots\subseteq O^{m_{t_{r}}}_{m_{t_{r}}}$$
   for all $r$. But $O^{m_{t_{r}}}_{m_{t_{r}}}\subseteq
   T_{m_{t_{r}}}$
   for all $r$, according to the definition of $\Upsilon$.
   It follows that $\cup_{n=1}^{\infty}O_{n}\subseteq\cup_{n=1}^{\infty}T_{n}.$
 $\QED$}
\end{description}

\begin{corollary}\label{cor:mgreals} There is a cofinal family ${\cal
   A}\subset\langle J_{{\reals}}\rangle$ such that TWO has a winning $2$-tactic
   in $MG({\cal A},J_{{\reals}})$. 
\end{corollary}

\begin{description}\item[Proof.]{Let ${\cal A}$ be the family of
   meager sets provided by Theorem \ref{important}.
   Thus, there is an order preserving function from $({\cal A},\subset)$
   to $(^{\omega}\omega,\ll)$. 
   But then  $({\cal
   A},\subset)\not\rightarrow(\omega-path)^{2}_{\omega/<\omega}$
holds, since  $(^{\omega}\omega,\ll)\not\rightarrow(\omega-path)^{2}_{\omega/<\omega}$
holds. By 
   Theorem \ref{important} the family ${\cal A}$ also satisfies the third
   hypothesis of Theorem \ref{mgth}. $\QED$}
\end{description}

\begin{corollary}\label{cor;localsmall} Let $J$ be a free ideal on an infinite
   set. If ${\cal A}$ is a family of sets in $\langle J\rangle$ such that: 
\begin{enumerate}
\item{${\cal A}$ is locally small,}
\item{for each $X\in{\cal A}$ there is a $Y\in{\cal A}$ such that $X\subset Y$,
   and}
\item{$({\cal A},\subset)$ is well-founded,}
\end{enumerate}
   then TWO has a winning $2$-tactic in $MG({\cal A},J)$.
\end{corollary}
\begin{description}\item[Proof.]{The proof is analogous to that of Corollary
  \ref{cor:mgreals}; now we refer to the proof of Theorem \ref{locsmallth}, 
  we observe that $\omega_{1}\leq{\frak b}$, and invoke Theorem \ref{mgth}.
  $\QED$} 
\end{description}

\begin{corollary}\label{cor:countablecof2} Let $\lambda\leq\kappa$ be
   infinite cardinal numbers such that:
\begin{enumerate}
\item{$\lambda$ has countable cofinality,}
\item{$\lambda^{+}\not\rightarrow(\omega-path)^{2}_{\omega/<\omega}$, and}
\item{$[\kappa]^{\leq\lambda}$ has the irredundancy property.}
\end{enumerate}
   Then there is a cofinal family ${\cal A}\subset[\kappa]^{\lambda}$ such
   that TWO has a winning 2-tactic in $MG({\cal A},[\kappa]^{<\lambda})$.
\end{corollary}

\begin{description}\item[Proof.]{Let ${\cal A}$ be a well-founded cofinal
   family in $[\kappa]^{\lambda}$ which is irredundant. Since there is a
   rank-function from ${\cal A}$ to $\lambda^{+}$ it follows from hypothesis
   {\em 2} that $({\cal A},\subset)\not\rightarrow(\omega-path)^{2}_{\omega/
   <\omega}$. From Corollary \ref{cor:countablecof} it follows that ${\cal A}$
   also
   satisfies the third hypothesis of Theorem \ref{mgth}. By that theorem TWO
   then has a winning $2$-tactic in the game $MG({\cal
   A},[\kappa]^{<\lambda})$. $\QED$} 
\end{description}

   Theorem \ref{slightprogress} shows
   that under certain circumstances there is for each $n$ a free ideal $J_n$
   and a cofinal family ${\cal A}_n\subset \langle J_n\rangle$ such that TWO
   does not have a winning $n$-tactic, but does have a winning $n+1$-tactic in
   $MG({\cal A}_n,J_n)$. We think that Theorem  \ref{slightprogress} indicates
   some relevance of the games as considered here for Telgarsky's Conjecture
   (see 3.4).

\begin{theorem}\label{slightprogress}Let $\lambda$ be an infinite cardinal
   number. If there is a linearly ordered set $(L,<)$ such that:
\begin{enumerate}
\item{$cof(L,<)>\omega$,}
\item{$(L,<)\rightarrow (\omega-path)^2_{\lambda/<\omega}$, but}
\item{$(L,<)\not\rightarrow(\omega-path)^3_{\lambda/<\omega}$,}
\end{enumerate} 
   then there is for each $n$ a free ideal $J_n$ and a cofinal family ${\cal
   A}_n\subset \langle J_n\rangle$ such that TWO does not have a winning
   $n$-tactic, but does have a winning $n+1$-tactic in $MG({\cal A}_n,J_n)$. 
\end{theorem}

\begin{description}\item[Proof.]{Let $\lambda$ and $(L,<)$ be as in the
   hypotheses. It follows from Propositions 3 and 4 of \cite{S2}
   that there is for each integer $m>1$ a linearly ordered set $(L_n,<_n)$ such
   that $\omega<cof(L_n,<_n)$ and:
\begin{equation}(L_n,<_n)\rightarrow(\omega-path)^{n}_{\lambda/<\omega}
\end{equation}
   but
\begin{equation}(L_n,<_n)\not\rightarrow(\omega-path)^{n+1}_{\lambda/<\omega}
\end{equation}
   Let $n>1$ and $(L_n,<_n)$ be fixed for the rest of the proof. We may
   assume that the underlying set, $L_n$, is disjoint from ${\EuScript
   P}({\EuScript P}(\lambda))\cup{\EuScript P}(\lambda)\cup\lambda$.\\

   Define a free ideal $J_n$ as follows: The underlying set on which $J_n$
   lives, say $S_n$, is $[\lambda]^{<\aleph_0}\cup L_n$. For each
   $\alpha\in\lambda$ let $X_{\alpha}$ be the set $\{Z\in
   [\lambda]^{<\aleph_0}: \alpha\not\in Z\}$.  Let ${\cal T}$ be
   $\{X_{\alpha}:\alpha\in\lambda\}$. Put a subset $X$ of
   $S_n$ in $J_n$ if:
\begin{quote} $X\cap[\lambda]^{<\aleph_0}$ is a subset of a union of finitely
   many elements of ${\cal T}$, and $X\cap L_n$ is bounded above.
\end{quote}
   
   Then the cofinality of $\langle J_n\rangle$ is $cof(L_n,<_n)$. Define
   ${\cal A}_n$ so that $X\in {\cal A}_n$ if:
   $$X\cap L_n=\{t\in L_n:t<z\} \mbox{for some $z\in L_n$.}$$
   Then ${\cal A}_n$ is cofinal in $\langle J_n\rangle$.\\
{\bf Claim 1:} TWO does not have a winning $n$-tactic in $MG({\cal
   A}_n,J_n)$.\\ 
   For let $\Phi$ be an $n$-tactic of TWO. For $x\in L_n$ put
   $V_x=[\lambda]^{<\aleph_0}\cup\{y\in L_n:y<_nx\}$. Define a partition
   $\Psi:[L_n]^{n}\rightarrow[\lambda]^{<\aleph_0}$ so that
$$ (\Phi(V_{x_1}) \cup \Phi(V_{x_1},V_{x_2})\cup \dots \cup
   \Phi(V_{x_1},\dots,V_{x_n}))\cap
   [\lambda]^{<\aleph_0}$$ is a  subset of $\cup\{X_{\alpha}:
   \alpha\in\Psi(\{x_{1},\dots,x_{n}\})\}$.  

   By $(1)$ we obtain an $\omega$-path $x_1<_nx_2<_n\ldots<_nx_k<_n\ldots$
   and a finite set $F\subset\lambda$ such that $\Psi(x_{j+1},\ldots,x_{j+n})
   \subseteq F$ for all $j$. For each $m$ we define:
   $O_m=[\lambda]^{<\aleph_0}\cup V_{x_m}$. Letting
   $(O_1,T_1,\ldots,O_k,T_k,\ldots)$ be the corresponding $\Phi$-play, we
   find that TWO has lost this play since
   $[\lambda]^{<\aleph_0}\cap(\cup_{m=1}^{\infty}T_m)\subseteq \cup_{\alpha\in
   F} X_{\alpha}\not=[\lambda]^{<\aleph_0}.$ 

   It follows that TWO does not have a winning $n$-tactic.\\

{\bf Claim 2:} TWO has a winning $n+1$-tactic in $MG({\cal
   A}_n,J_n)$.\\ 

   First observe that $\cup_{\alpha\in F}X_{\alpha}=[\lambda]^{<\aleph_0}$ 
   whenever $F$ is an infinite subset of $\lambda$. 

   Here is a definition of an $n+1$-tactic for TWO in this game: Let
   $\{t_{\alpha}:\alpha<\lambda\}$ enumerate $[\lambda]^{<\aleph_0}$
   bijectively. Let $\Phi:[L_n]^{n+1}\rightarrow\lambda$ be a coloring which
   witnesses that
   $(L_n,<_n)\not\rightarrow(\omega-path)^{n+1}_{\lambda/<\omega}$. For each
   $X$ in ${\cal A}_n$ let $\phi_X$ be that element of $L_n$ for which
   $X\cap L_n=\{t\in L_n:t<\phi_X\}$. 

   For $U_1\subset\ldots\subset U_{n+1}$ elements of ${\cal A}_n$, observe that
   $\phi_{U_1}\leq\ldots\leq\phi_{U_{n+1}}$. For $X\subset Y$ sets in ${\cal
   A}_n$ such that $X\cap[\lambda]^{<\aleph_0}\neq Y\cap[\lambda]^{<\aleph_0}$
   we set $\Psi(X,Y)=\min\{\alpha:t_{\alpha}\in Y\backslash X\}$.

   Let $U_1\subset\ldots\subset U_{n+1}\in {\cal A}_n$ be given. We define:\\
\begin{enumerate}
\item{$G(U_1,\ldots,U_j)=\emptyset$ when $j<n+1$,}
\item{$G(U_1,\ldots,U_{n+1})=X_{\alpha}\cup (L_n\cap U_{n+1})$ when
   $\phi_{U_1}<\ldots<\phi_{U_{n+1}}$, and
   $\Phi(\{\phi_{U_1},\ldots,\phi_{U_{n+1}}\}) = \alpha$,}
\item{$G(U_1,\ldots,U_{n+1})=X_{\alpha}\cup (L_n\cap U_{n+1})$ where $\alpha$
   is minimal such that $t_{\alpha}\in U_{i+1}\backslash U_i$ for some
$i\leq n$,
   otherwise.} 
\end{enumerate}

   We show that $G$ is a winning $n+1$-tactic for TWO. Thus, let
$$(O_1,T_1,\ldots, O_m,T_m,\ldots)$$ 
   be a $G$-play of the game. For typographical convenience we define:
\begin{enumerate}
\item{$x_i=\phi_{O_i}$ for each $i$, and}
\item{$\alpha_i=\Psi(O_i,O_{i+1})$ for each $i$ for which this is defined.}
\end{enumerate}
   There are two cases to consider.\\
CASE 1: $\{i:x_i=x_{i+1}\}$ is finite.\\
  Choose $m$ such that $x_i<x_{i+1}$ for all $i\geq m$. Then the set
  $$\{\Phi(\{x_{m+k+1}, \ldots, x_{m+k+n+1}\}): k=1, 2, \ldots\}$$ is
  an infinite
   subset of $\lambda$ and it follows from 2. in the definition of $G$ that this
   play is won by TWO.\\ 
CASE 2: $\{i:x_i=x_{i+1}\}$ is infinite.\\
  Then the set $\{i:\Psi(O_i,O_{i+1})\mbox{ is defined}\}$ is infinite. But
  then it follows from 3. in the definition of $G$ that TWO wins this play.
$\QED$} \end{description}

   The hypotheses of Theorem \ref{slightprogress} are realized under
   any of the following axiomatic circumstances (one uses Corollary 27 and
   Proposition 29 of \cite{S2}):

\begin{enumerate}
\item{$2^{<{\frak c}}={\frak c}+EH$,}
\item{${\frak c}<2^{\aleph_1}$, i.e., the negation of $LH$ (Lusin's
   second Continuum Hypothesis)}
\item{There is an infinite regular cardinal number $\kappa$ such that
   $2^{\kappa}=\kappa^+$.}
\end{enumerate}  

   For the case when $\lambda=\omega$, the
   example constructed in the proof of Theorem \ref{slightprogress} shows
   that hypothesis 2 of Theorem \ref{mgth} is to some extent necessary.
   This is because:
\begin{enumerate}
\item{${\cal A_n}$ has the coherent decomposition 
   property: For choose $\alpha_1<\alpha_2<\dots<\alpha_n<\dots$ from
   $\omega$, and set $T_m=X_{\alpha_1}\cup\dots\cup X_{\alpha_m}$ for
   each $m$. Then
   $[\omega]^{<\aleph_0}=\cup_{m=1}^{\infty}X_{\alpha_m}$, and
   $X_{\alpha_j}\subseteq X_{\alpha_i}$ for $j<i$.
   For $A\in {\cal A_n}$ we put $A_m=(A\cap T_m)\cup(A\cap L_n)$.}
\item{$({\cal
   A}_m,\subset)\rightarrow(\omega-path)^m_{\omega/<\omega}$, but}
\item{$({\cal
   A}_m,\subset)\not\rightarrow(\omega-path)^{m+1}_{\omega/<\omega}$.}
\end{enumerate}

   At this point it is an open problem whether the hypotheses (and for that
   matter the conclusion) of Theorem \ref{slightprogress} are
   satisfied simply in the theory ZFC (see Problem 9 of \cite{S2}).

\subsection{The game MG(J)}

   $MG(J)$ denotes the version of $MG({\cal A},J)$ where $\langle
   J\rangle = {\cal A}$. In Problem 1 of
   \cite{S1} it was asked whether there is for each $k$ a free ideal
   $J_k$ such that TWO does not have a winning $k$-tactic in $MG(J_k)$,
   but does have a winning $k+1$-tactic in $MG(J_k)$. This problem is
   still open. In \cite{S1}, Corollary 10, it was proven that TWO does
   not have a winning 
   $2$-tactic in the game $MG(J_{\reals})$, but that TWO has a winning
   $3$-tactic in $MG(J_{\reals})$ if for example the Continuum
   Hypothesis is assumed. We now extend these results in two
   directions. 
\begin{enumerate}
\item{In Problem 
   3 of that paper it was asked if player TWO has a winning $3$-tactic if
   instead of the Continuum Hypothesis one uses the theory $ZFC+MA+EH+\neg CH$,
   which is explained below. We now show that the answer is
   affirmative.}
\item{We identify circumstances under which TWO does not have a
   winning $k$-tactic in $MG(J_{\reals})$ for any $k$; combining this
   with a consistency result of Todorcevic (given in the appendix), it
   follows that it is also consistent that there is no $k$ for which TWO
   has a winning $k$-tactic in $MG(J_{\reals})$.}
\end{enumerate}

   It follows that the existence of a winning $k$-tactic for TWO in
   $MG(J_{\reals})$ is not decided by the axioms of traditional set theory.
   One might now wonder if it is consistent that for example TWO does
   not have a winning $3$-tactic in $MG(J_{\reals})$, but does have a
   winning $4$-tactic? This is not possible since a theorem of \cite{S3}
   implies that either TWO has a winning $3$-tactic, or else there is no
   $k$ such that TWO has a winning $k$-tactic in $MG(J_{\reals})$.

   Let $EH$ (which abbreviates {\em Embedding Hypothesis}) denote the
   statement: 	
$$\mbox{{\em every linearly ordered set of cardinality $\leq {\frak c}$
   embeds in }$(^{\omega}\omega,\ll)$.}
$$ 
   The hypothesis $EH$ is a consequence of the Continuum Hypothesis.
   Laver has proven (\cite{L}) that the theory $ZFC+EH+\neg CH$ is
   consistent, and Woodin (\cite{W}, pp. 31-47), extending this, has
   proven the consistency of the theory $ZFC+MA+EH+\neg CH$. This
   theory implies that $2^{<\frak c}={\frak c}+EH$, which in turn is
   strong enough to prove that the partition relation $$({\cal
   P}({\frak c}),\subset)\not\rightarrow(\omega-\mbox{path})^{3}_{\omega
   /<\omega}$$ holds (see \cite{S2}, top of p. 60). Thus we have:

\begin{prop}\label{problem3} The theory $``ZFC+\neg CH+$ TWO has a
   winning 3-tactic in $MG(J_{{\reals}})$" is consistent.
\end{prop}
\begin{description}\item[Proof.]{Consider any model of $ZFC+EH+\neg
   CH+2^{<\frak c}={\frak c}$ in which $\langle J_{\reals}\rangle$ has
   a cofinal chain. Let ${\EuScript C}$ denote this cofinal chain. By
   Theorem \ref{important} we may assume that this cofinal chain has a
   coherent decomposition and that it satisfies the partition relation
$({\EuScript C},\subset)\not\rightarrow(\omega-path)^2_{\omega/<\omega}.$ 
   Since we also have 
$({\cal
   P}({\frak c}),\subset)\not\rightarrow(\omega-\mbox{path})^{3}_{\omega
   /<\omega}$ it follows that:
\begin{enumerate}
\item{$(\langle J_{\reals}\rangle,\subset)\not \rightarrow
   (\omega-path)^3_{\omega/<\omega}$, and}
\item{The family $\langle J_{\reals}\rangle$ has a coherent
   decomposition.}
\end{enumerate}

   Theorem \ref{mgth} implies that TWO has a winning $3$-tactic in
   $MG(\langle J_{\reals}\rangle, J_{\reals})$. This completes the proof
   of the proposition. $\QED$}
\end{description}

   Indeed, our proof of Proposition \label{problem3} shows more
   generally that if $J$ is a free ideal on a set of cardinality at most
   ${\frak c}$, and if $\JA$ has a cofinal chain and the coherent
   decomposition property, then the theory $ZFC+EH+2^{<{\frak c}}={\frak
   c}$ proves that TWO has a winning $3$-tactic in $MG(J)$ . This
   generalizes Theorem 8(a) of \cite{S1}.

   Next we give hypotheses under which there is no $k$ for which TWO has a
   winning $k$-tactic in $MG(J_{{\reals}})$. In the appendix we give a
   proof that these hypotheses are consistent with $ZFC$. This consistency
   result is due to Todorcevic.

\begin{theorem}\label{noktactic} Assume that
   $cof(J_{{\reals}},\subset)=\lambda$ and that the partition relation
   $({\EuScript
P}({\frak c}),\subset)\rightarrow(\omega-path)^{3}_{\lambda/<\omega}$
   holds. Then there is no $k$ for which TWO has a winning $k$-tactic in
   $MG(J_{\reals})$. 
\end{theorem}

\begin{description}\item[Proof.]{Let $k$ as well as a $k$-tactic $F$ for $TWO$
   be given. Let $X$ be a nowhere dense subset of cardinality ${\frak c}$
   of ${\reals}\backslash{\rationals}$. Let ${\cal
   A}=\{A_{\alpha}:\alpha<\lambda\}$ be a bijectively enumerated cofinal
   subfamily of $J_{\reals}$.

   Define a partition $\Phi:[{\EuScript
   P}(X)]^{k}\rightarrow\lambda$ so that 
   $$\Phi(\{X_{1},\cdots,X_{k}\})=\beta$$ where $\beta$ is minimal
   such that $$F({\rationals}\cup X_{1})\cup\cdots\cup
   F({\rationals}\cup X_{1},\cdots,{\rationals}\cup X_{k})\subset
   A_{\beta}.$$

   Since $({\EuScript
   P}({\frak c}),\subset)\rightarrow(\omega-path)^{3}_{\lambda/<\omega}$,    it follows that
   $({\EuScript
   P}({\frak c}),\subset)\rightarrow(\omega-path)^{k}_{\lambda/<\omega}$ (   see 
   \cite{S2}, Proposition 36). Accordingly, choose a finite set $G\subset
   \lambda$ and an increasing $\omega$-sequence
   $X_{1}\subset X_{2}\subset\cdots$ of subsets of  $X$ such that
   $\Phi(\{X_{j+1},\cdots,X_{j+k}\})\in
   G$ for all $j$. Put $O_{n}=X_{n}\cup{\rationals}$ for all $n$. Let $B$
   be the nowhere dense set $\cup\{A_{\alpha}:\alpha\in G\}$. Also define
   $T_{j}=F(O_{1}\cdots,O_{j})$ for $j\leq k$, and
   $T_{j+k}=F(O_{j+1},\cdots,O_{j+k})$ for all $j$. Then
   $$(O_{1},T_{1},O_{2},T_{2},\cdots)$$ is an $F$-play of $MG(J_{\reals})$ for
   which ${\rationals}\subset\cup_{n=1}^{\infty}O_{n}$ and
   $\cup_{n=1}^{\infty}T_{n}\subseteq B$. Since $B$ is nowhere dense,
   ${\rationals }\backslash B\neq \emptyset$. It follows that TWO has lost this
   play. $\QED$} 
\end{description}

   We now consider games of the form $MG([\kappa]^{<\lambda})$. In
   Proposition 15 of \cite{S1} it was shown that if TWO has a winning
   $k$-tactic in this game for some $k$, then TWO in fact has a winning
   $3$-tactic. It is not known if $``3"$ is optimal (this is Problem
   $7$ of \cite{S1}). It also follows from \cite{S1}, Proposition $5$, that if
   $\lambda\rightarrow(\omega-path)^{2}_{\omega/<\omega}$, then TWO does not
   have a winning $k$-tactic in this game for any $k$. We now present
   slightly sharper results.

\begin{theorem}\label{th:countablecof} Let $\lambda$ be an uncountable cardinal
   number of countable cofinality. Let $k>1$ be an integer. The following
   statements are equivalent: 
\begin{enumerate}
\item{Player TWO has a winning $k$-tactic in the game
   $MG([\lambda^{+}]^{<\lambda})$.}
\item{$([\lambda^+]^{\leq\lambda},\subset)\not\rightarrow
   (\omega-path)^k_{\omega/<\omega}$.}
\item{$\lambda^{+}\not\rightarrow(\omega-path)^{2}_{\omega/<\omega}$ and
   $({\cal
   P}(\lambda),\subset)\not\rightarrow(\omega-path)^{k}_{\omega/<\omega}$.}
\end{enumerate}
\end{theorem}
\begin{description}\item[Proof.]{By Theorem 1 and Proposition 15 of \cite{S1}
   we may assume that $k\in\{2,3\}$. Let $\lambda_1<\ldots<\lambda_n<\ldots$ be
   a sequence of cardinal numbers converging to $\lambda$.\\

$1.\Rightarrow 2.$\\ Let $F$ be a winning $k$-tactic for TWO in
   $MG([\lambda^+]^{<\lambda})$. Put $S=\lambda^+\backslash \lambda$.
   Define a coloring $\Phi:[[S]^{\leq\lambda}]^k\rightarrow\omega$
   so that 
$$\Phi(X_1,\dots,X_k)=\min\{n:|F(\lambda\cup X_1,\dots,\lambda\cup
   X_k|\leq \lambda_n\}.
$$
   Since $F$ is a winning $k$-tactic for TWO, $\Phi$ is a coloring
   which witnesses the partition relation in $2$.

$2.\Rightarrow 1.$\\

   According to Corollary \ref{cor:countablecof},
$([\lambda^+]^{<\lambda}, \subset)$ has the coherent decomposition
property. Since $[\lambda^+]^{\leq\lambda}$ has a cofinal chain it
follows that this family of sets itself has a coherent decomposition.
The partition property in $2$ implies that the family
$[\lambda^+]^{\leq\lambda}$ satisfies the hypotheses of Theorem
\ref{mgth}; thus TWO has a winning $k$-tactic in
$MG([\lambda^+]^{<\lambda})$.

   The equivalence of 2. and 3. is also easy to establish. $\QED$}
\end{description}

\begin{corollary}Let $\lambda$ be an uncountable cardinal number of countable
   cofinality. Assume
   $ZFC+EH+\lambda<{\frak c}+{\frak c}=2^{<{\frak c}}$. Then TWO has
   a winning $2$-tactic in $MG([\lambda^{+}]^{<\lambda})$.
\end{corollary}
\begin{description}\item[Proof.]{The hypothesis
   $EH+{\frak c}=2^{<{\frak c}}$ implies that both $\lambda^{+}$ and
   $({\EuScript P}(\lambda),\subset)$ embed in $(^{\omega}\omega,\ll)$ for any
   $\lambda<{\frak c}$. It then follows from Corollary 13 of
   \cite{S2} that the partition relations in $3.$ of Theorem 
   \ref{th:countablecof} hold for $k=2$ for each $\lambda<{\frak c}$. $\QED$}
\end{description}

\subsection{The game SMG(J)}

   For a free ideal $J$ on an infinite set $S$, the game $SMG(J)$ (read
   ``strongly monotonic game on J") is defined so that an $\omega$-sequence
   $(O_{1},T_{1},\dots,O_{n},T_{n},\dots)$ is a play if for each $n$,
\begin{enumerate}
\item{$O_{n}\in\langle J\rangle$ is player ONE's move in inning $n$,}
\item{$T_{n}\in J$ is player TWO's move in inning $n$, and}
\item{$O_{n}\cup T_{n}\subseteq O_{n+1}$.}
\end{enumerate}
   Player TWO wins this play if
   $\cup_{n=1}^{\infty}O_{n}=\cup_{n=1}^{\infty}T_{n}.$

   Throughout this section we assume that $\langle J\rangle$ is a
   proper ideal on $S$.

\begin{theorem}\label{smgth} Let $J\subset{\EuScript P}(S)$ be a free
ideal and let 
${\cal A}$ be a cofinal subfamily of $\langle J\rangle$ such that:
\begin{enumerate}
\item{TWO has a winning $k$-tactic in $MG({\cal A},J)$,}
\item{there are functions $\Phi_1:\langle J\rangle\rightarrow J$ and
$\Phi_2:\langle J\rangle\rightarrow{\cal A}$ such that:
   \begin{enumerate}
   \item{$A\subset \Phi_2(A)$ for each $A\in \langle J\rangle$, and}
   \item{$\Phi_2(A)\subset\Phi_2(B)$ whenever $A\cup\Phi_1(A)\subseteq
   B\in \langle J\rangle$.}
   \end{enumerate}}
\end{enumerate}
   Then TWO has a winning $2$-tactic in $SMG(J)$.
\end{theorem}

\begin{description}\item[Proof]{Let ${\cal A}$, $\Phi_1$ and $\Phi_2$
   be as in the hypotheses. For each $A$ in $\langle J\rangle$ define
   $(A_1,\dots,A_k)$ so that $A_1=\Phi_2(A)$ and $A_{j+1}=\Phi_2(A_j)$
   for each $j<k$. Also define:
   $\Psi(A)=\Phi_1(A)\cup\Phi_1(A_1)\cup\dots\cup \Phi_1(A_k)$. 

   Let $F$ be a winning $k$-tactic for TWO in
   $MG({\cal A},J)$. Define a $k$-tactic, $G$, for TWO as follows. Let
   $A\subset B$ be given.\\

CASE 1: $G(A)=F(A_1)\cup\dots\cup F(A_1,\dots,A_k)\cup\Psi(A)$.\\

CASE 2: If $A_k\subset B_1$, we let $G(A,B)$ be the set
$$F(A_2,\dots,A_k,B_1)\cup
F(A_3,\dots,A_k,B_1,B_2)\cup\dots\cup F(B_1,\dots,B_k)\cup\Psi_1(B).$$

CASE 3: Otherwise we put $G(A,B)=G(B).$\\

   Then $G$ is a winning $2$-tactic for TWO in $SMG(J)$. For let 
$$(O_{1},T_{1},\dots,O_{n},T_{n},\dots)$$
   be a play of $SMG(J)$ during which TWO followed the $2$-tactic $G$.
   For each $j$ we put
   $M_j^1=\Phi_2(O_j),\dots,M_j^k=\Phi_2(M_j^{k-1})$. An inductive
   computation shows that

\begin{itemize} 
\item{$(M_1^1,M_1^2,\dots,M_1^k,M_2^1,M_2^2,\dots,M_2^k,\dots)$ is a
   sequence of legal moves for $ONE$ in 
   the game $MG({\cal A},J)$, and that}
\item{\begin{enumerate}
\item{$F(M_1^1)\cup\dots\cup F(M_1^1,\dots,M_1^k)\subseteq T_1$, and} 
\item{$F(M_j^1,\dots,M_j^k)\cup F(M_j^2,\dots,M_j^k,M_{j+1}^1)\cup \dots
   \cup F(M_j^k,M_{j+1}^1,\dots,M_{j+1}^{k-1})\\ \subseteq T_{j+1}$
   for each $j$.}
\end{enumerate}}
\end{itemize} 

   Since $F$ is a winning $k$-tactic for TWO in the game
   $MG({\cal A},J)$, and since $\cup_{n=1}^{\infty}O_n\subseteq
   \cup_{n=1}^{\infty}M_n^1$, TWO won the given play of
   $SMG(J)$. $\QED$}
\end{description}

   The next corollary solves Problems 10 and 11 of \cite{S1}. 

\begin{corollary}\label{corsmg} Player $TWO$ has a winning 2-tactic in the
   game $SMG(J_{{\reals}})$.
\end{corollary}

\begin{description}\item[Proof.]{ Fix, by Corollary \ref{cor:mgreals},
   a cofinal family ${\cal A}\subset \langle J_{\reals}\rangle$ such that
   TWO has a winning $2$-tactic in $MG({\cal A},J_{\reals})$. 

   We define 
$\Phi_1:\langle J_{\reals}\rangle\rightarrow J_{\reals}$
and
$\Phi_2:\langle J_{\reals}\rangle\rightarrow {\cal A}$
as follows:\\

Fix $X\in\langle J_{\reals}\rangle$, and choose a sequence
$(X_0,X_1,\dots, X_n,\dots)$ such that:
\begin{enumerate}
\item{$X_0=X$,}
\item{$X_{n+1}\in{\cal A}$ and $\naturals\cdot X_n\subseteq X_{n+1}$}
\end{enumerate}
   for each $n$. Put $\Phi_2(X)=\cup_{n=1}^{\infty}X_n.$

   Fix $X\in\langle J_{\reals}\rangle$ and let $\Phi_1(X)$ be a nowhere
   dense set for which $\Phi_2(X)\subset \naturals\cdot\Phi_1(X)$.

   Then ${\cal A}$, $\Phi_1$ and $\Phi_2$ are as required by Theorem
   \ref{smgth}. $\QED$}
\end{description}

\begin{corollary}\label{Th21} For each of the ideals $J_n$ constructed
in the proof of Theorem \ref{slightprogress}, TWO has a winning
$2$-tactic in $SMG(J_n)$.
\end{corollary}

\begin{description}\item[Proof.]{Let ${\cal A}_n$ be as in the proof
   of Theorem \ref{slightprogress}. For each $X\in\langle J_n\rangle$ we
   let $\Phi_2(X)$ be an element of ${\cal A}_n$ which contains it, and
   we let $\Phi_1(X)=\{a_X\}$ where $a_X\in L_n\backslash\Phi_2(X)$. Then
   ${\cal A}_n$, $\Phi_1$ and $\Phi_2$ are as required by Theorem
   \ref{smgth}. $\QED$}
\end{description}

   Before giving another application of Theorem \ref{smgth} we give an
   example of free ideals $J$ which show that TWO does not always
   have a winning $k$-tactic in the game $SMG(J)$ for some $k$. These
   examples are also relevant to the material of the next section. The symbol
   $M(\omega,2)$ denotes the smallest ordinal $\alpha$ for which the partition
   relation $\alpha\rightarrow(\omega-path)^{2}_{\omega/<\omega}$ holds.
   $M(\omega,2)$ is a regular cardinal less than or equal to
   ${\frak c}^{+}$. It in fact satisfies the partition relation
   $M(\omega,2)\rightarrow(\omega-path)^{n}_{\omega/<\omega}$ for all $n$. Let
   $\kappa$ be an initial ordinal number. It is consistent that $M(\omega,2)$
   is equal to $\aleph_{2}$ while ${\frak c}$ is larger than
   $\kappa$ (this is yet another result of Todorcevic). 

\begin{theorem}\label{counterexample}Let $\lambda$ be a cardinal number of
   countable cofinality and let $\kappa$ be a cardinal number larger
   than $\lambda$. If $M(\omega,2) \leq \lambda^+$, then
   there is no $k$ such that player TWO has a winning $k$-tactic in
   $SMG([\kappa]^{<\lambda})$.
\end{theorem}

\begin{description}\item[Proof.]{Let $F$ be a $k$-tactic for TWO. 

   Player ONE's counter-strategy will be to play judiciously chosen subsets
   from $\kappa$. We first single out those
   sets from which ONE will make moves.  

   Choose sets $S_0\subset S_1\subset\dots\subset S_{\alpha}\subset
\dots\in[\kappa]^{\lambda}$ for $\alpha<\lambda^+$ such that:
\begin{enumerate}
\item{$\lambda\subset S_0$,}
\item{$\cup\{F(S_{i_{1}},\dots,S_{i_j}):j\leq k,\
   i_1<\dots<i_j<\alpha\}\subset S_{\alpha}$ for each
   $0<\alpha<\lambda^+$.}
\end{enumerate}

   Now let $\lambda_{1}<\lambda_{2}<\dots<\lambda$ be an increasing sequence of
   regular cardinal numbers converging to $\lambda$. Define a function
$
\Gamma:[\lambda^+]^{k}\rightarrow\omega
$
     so that 
$$  \Gamma(\xi_{1},\dots,\xi_{k})=\min\{m:|F(S_{\xi_{1}},\dots,S_{\xi_{k}})|
    \leq\lambda_{m}\}.
$$

   Then, on account of the relation $M(\omega,2)\leq\lambda^{+}$, choose an $m<
   \omega$ and a sequence $\alpha_{k+1}<\dots<\alpha_{k+m}<\dots$ from
   $\lambda^+$
   such that $\Gamma(\alpha_{j+1},\dots,\alpha_{j+k})\leq m$ for all $j$.

   Consider the sequence 
$$(S_{\alpha_{1}}, F(S_{\alpha_{1}}),\dots,S_{\alpha_{k}},
   F(S_{\alpha_{1}},\dots,S_{\alpha_{k}}),\dots,S_{\alpha_{k+m}},
   F(S_{\alpha_{1+m}},\dots,S_{\alpha_{k+m}}),\dots).
$$
   It is a play of the game $SMG([\kappa]^{<\lambda})$ during which TWO
   used the $k$-tactic $F$. To see that TWO lost this play, let $T_{k}$
   denote TWO's $k$-th move. The choice of the sequence $\alpha_{k+1}<\dots$
   implies that $\cup_{n=1}^{\infty}T_{n}$ has cardinality less than $\lambda$.
   The union of the sets played by ONE has cardinality $\lambda$; TWO didn't
   catch up with ONE. $\QED$}
\end{description}

\begin{corollary}\label{cor:countablecof4} For $\omega=cof(\lambda)\leq\lambda<\kappa$ cardinal numbers with
$cof([\kappa]^{\leq\lambda},\subset)=\kappa$, the following
statements are equivalent:
\begin{enumerate}
\item{TWO has a winning $2$-tactic in $SMG([\kappa]^{<\lambda})$.}
\item{$\lambda^+\not\rightarrow(\omega-path)^2_{\omega/<\omega}$.}
\end{enumerate}
\end{corollary}

\begin{description}\item[Proof.]{It follows from Theorem
   \ref{counterexample} that $1.$ implies $2.$ \\

   That $2.$ implies $1.$:\\

   By the cofinality hypothesis and by $2.$ we find, according to
   Corollary \ref{cor:countablecof2}, a well-founded cofinal family
   ${\cal A}$ such that TWO has a winning $2$-tactic in $MG({\cal
A},[\kappa]^{<\lambda})$. We may assume that there is an enumeration
$\{A_{\alpha}:\alpha<\kappa\}$ of ${\cal A}$ for which $\alpha\in
A_{\alpha}$ for each $\alpha$. Define $\Phi_1$ and $\Phi_2$ as follows:\\

For $X\in[\kappa]^{\leq\lambda}$ define a sequence
$(X_0,\dots,X_m,\dots)$ such that:

\begin{enumerate}
\item{$X_0=X$, and}
\item{$X_{n+1}=\cup_{\alpha\in X_n}A_{\alpha}$}
\end{enumerate}
   for each $n$.

Choose $\Phi_2(X)\in{\cal A}$ such that
$\cup_{n<\omega}X_n\subseteq\Phi_2(X)$.\\

Pick $z_X\in(\kappa\backslash\Phi_2(X))$ and pick $\rho_X$ minimal
such that $\rho_X\not\in\Phi_2(X)$, and $\Phi_2(X)\subset A_{\rho_X}$. Put
$\Phi_1(X)=\{z_X, \rho_X\}$.\\

Then ${\cal A}$, $\Phi_1$ and $\Phi_2$ are as required by Theorem
\ref{smgth}. $\QED$}
\end{description}

This result will be discussed at greater length after Theorem \ref{3tacticth}. 

   We finally mention that it is still unknown whether there is for
   each $m$ a free ideal $J_m$ such that TWO does not have a winning
   $m$-tactic, but does have a winning $m+1$-tactic in $SMG(J_m)$.
   This is Problem 9 of \cite{S1}. In this connection it is worth
   noting the following 
   relationship between winning $k$-tactics in $MG(J)$ and winning
   $m$-tactics in $SMG(J)$. The proof uses ideas as in the proof of
   Theorem \ref{smgth}.

\begin{theorem} If TWO has a winning $k$-tactic in $MG(J)$, then TWO
   has a winning $2$-tactic in $SMG(J)$.
\end{theorem}

\subsection{The game VSG(J)}

   For a free ideal $J$ on an infinite set $S$, the game $VSG(J)$ (read
   `` very strong game on J") is defined so that an $\omega$-sequence
   $(O_{1},(T_{1},S_{1}),\dots,O_{n},(T_{n},S_{n}),\dots)$ is a play if for
   each $n$, 
\begin{enumerate}
\item{$O_{n}\in\langle J\rangle$ is player ONE's move in inning $n$,}
\item{$(T_{n},S_{n})\in J\times\langle J\rangle$ is player TWO's move in inning
   $n$, and} 
\item{$O_{n}\cup T_{n}\cup S_{n}\subseteq O_{n+1}$.}
\end{enumerate}
   Player TWO wins this play if
   $\cup_{n=1}^{\infty}O_{n}=\cup_{n=1}^{\infty}T_{n}.$

   We assume for this section that $\langle J\rangle$ is also a
   proper ideal on $S$. Given a cofinal family ${\cal A}\subset\langle
   J\rangle$, we may assume whenever convenient that ONE is playing from
   ${\cal A}$ in the game $VSG(J)$. It is clear that if TWO has a winning
   $k$-tactic in $SMG(J)$, then TWO has a winning $k$-tactic is $VSG(J)$.
   The converse is not so clear.

\begin{problem} Let $J$ be a free ideal on a set $S$ and let $k$ be a
   positive integer. Is it true that if TWO has a winning $k$-tactic in
   $VSG(J)$, then TWO has a winning $k$-tactic in $SMG(J)$?
\end{problem}

   In the next theorem we find a partial converse.

\begin{theorem}\label{SMG,VSG} Let $J$ be a free ideal on a set $S$ and let $k$
   be a positive integer. If $add(\langle J\rangle,\subset) = cof(\langle
   J\rangle,\subset)$, then the following statements are equivalent:
\begin{enumerate}
\item{TWO has a winning $2$-tactic in $SMG(J)$.}
\item{TWO has a winning $k$-tactic in $SMG(J)$.}
\item{TWO has a winning $k$-tactic in $VSG(J)$.}
\end{enumerate}
\end{theorem}
\begin{description}\item[Proof]{That $1.$ and $2.$ are equivalent: This is
Theorem 19 of \cite{S1}.\\
That $2.$ implies $3.$: Let $F$ be a winning
   $k$-tactic for TWO in $SMG(J)$. Define $G$ so that 
$$G(A_1,\dots,A_j)=(F(A_1,\dots,A_j),A_j\cup F(A_1,\dots,A_j))$$
   for $j\leq k$. Then $G$ is a winning $k$-tactic for TWO in $VSG(J)$.\\
That $3.$ implies $2.$: Let $G$ be a winning $k$-tactic for TWO in $VSG(J)$.
Then choose a sequence $(M_\xi:\xi<cof(\langle J\rangle,\subset))$ such that:
\begin{enumerate}
\item{$M_\xi\subset M_\nu$ for $\xi<\nu<cof(\langle J\rangle,\subset)$ and}
\item{$\{M_\xi:\xi<cof(\langle J\rangle,\subset)\}$ is cofinal in $\langle
  J\rangle$.}
\end{enumerate}

   Now $cof(\langle J\rangle,\subset)$ is a regular uncountable cardinal
   number. We may thus further assume that the sequence $(M_\xi:\xi<cof(\langle
   J\rangle,\subset))$ has been chosen such that if $(U,T)=G(M_{\xi_1},\dots,
   M_{\xi_j})$, then $U\cup T\subset M_{\eta}$ for all $\xi_j<\eta< cof(\langle
   J\rangle,\subset)$.

   For each $X\in \langle J\rangle$ define $\alpha(X)=\min\{\xi:X\subset
   M_\xi\}$. For each $\xi$ choose $z_\xi\in S\backslash M_\xi$. We now
   define a $k$-tactic, $F$, for TWO in $SMG(J)$. 

   Let $X_1\subset\dots\subset X_j\in\langle J\rangle$ for a $j\leq k$ be
   given. \\ 
CASE 1:$\alpha(X_1)<\dots<\alpha(X_j)$. Let $(U,T)=G(M_{\alpha(X_1)},\dots,
   M_{\alpha(X_j)})$ and define $F(X_1,\dots,
   X_j)=U\cup\{z_{\alpha(X_j)+1}\}$.\\
CASE 2: Otherwise, set $F(X_1,\dots, X_j)=\{z_{\alpha(X_j)+1}\}$.
   Then $F$ is a winning $k$-tactic for TWO in $SMG(J)$. $\QED$}
\end{description}

There is the following analogue of Theorem \ref{smgth} for the very
strong game:

\begin{prop}\label{mgtovsg}Let $J$ be a free ideal on a set $S$. If there is a
   cofinal family ${\cal A}\subset\langle J\rangle$ such that TWO has a winning
   $k$-tactic in $MG({\cal A},J)$, then TWO has a winning $2$-tactic in
$VSG(J)$. \end{prop}

\begin{description}\item[Proof.]{Let ${\cal A}\subset\langle J\rangle$ be a
   cofinal family such that TWO has a winning $k$-tactic in $MG({\cal A},J)$.
   We will define a winning 2-tactic for TWO for the game $VSG(J)$. To this
   end, choose a winning $k$-tactic, $F$, for TWO for the game
$MG({\cal A},J)$. 
   For each $X\in\langle J\rangle$ choose a set
   $A_1(X)\subset\dots\subset A_k(X)$ from ${\cal A}$
   such that $X\subset A_1(X)$, and choose $\Psi(X)$ from ${\cal A}$
   such that $A_k(X)\subset\Psi(X)$.

   Let $X\subset Y$ be sets from $\langle J\rangle$.\\
CASE 1: $G(X)=(F(A_1(X))\cup\dots\cup
F(A_1(X),\dots,A_k(X)),\Psi(X))$.\\
CASE 2: Define $G(X,Y)$ so that:
\begin{enumerate}  
\item{ $G(X,Y)= (F(A_2(X),\dots,A_k(X),A_1(Y) \cup\dots\cup F(A_1(Y),\dots,A_k(Y)), \Psi(Y))$ if
   $\Psi(X)\subset Y$, and} 
\item{$G(X,Y)= G(Y)$ otherwise.}
\end{enumerate}

   Then $G$ is a winning $2$-tactic for TWO in $VSG(J)$. For let 
   $(O_{1},(T_{1},S_{1}),O_{2},(T_{2},S_{2}),\dots)$ be a play of $VSG(J)$
   such that $(T_1,S_1)=G(O_1)$ and
   $(T_{n+1},S_{n+1})=G(O_n,O_{n+1})$ for all $n$. Then $S_n=\Psi(O_n)$
   and $A_k(O_n)\subset A_1(O_{n+1})$ for each $n$. An inductive
   computation, using this information, shows that TWO won this play
   of $VSG(J)$. $\QED$} 
\end{description}

   Combining Theorem \ref{SMG,VSG} and Theorem \ref{counterexample} we see
   that TWO does not always have a winning $k$-tactic in games of the form
   $VSG(J)$. Combining Theorem \ref{SMG,VSG} and Corollary
   \ref{cor:countablecof4} we obtain another game-theoretic characterization of
   the
   partition relation
   $\lambda^{+}\rightarrow(\omega-path)^{2}_{\omega/<\omega}$ when $\lambda$ is
   an uncountable cardinal of countable cofinality.

   Analogous to the case of the ideal of countable subsets of an infinite set, 
   there is for each uncountable cardinal number $\lambda$ which is of
   countable cofinality, a proper class of cardinals $\kappa$ for which the
   ideal $[\kappa]^{\leq\lambda}$ has the irredundancy property. It is also a
   consequence of $MA+{\frak c}>\lambda$ that the partition relation
   $\lambda^{+} \not \rightarrow(\omega-path)^{2}_{\omega/<\omega}$ holds.
   Accordingly it is consistent that there is a proper class of cardinals
   $\kappa$ such that TWO has a winning $2$-tactic in the game
   $VSG([\kappa]^{\leq\lambda})$. The following problem (to be compared with
   the upcoming Conjecture \ref{conj:countable-finite}) is open.

\begin{problem} Let $\lambda$ be an uncountable cardinal of countable
   cofinality. Is it true that if TWO has a winning $2$-tactic in the game
   $VSG([\lambda^{+}]^{<\lambda})$, then TWO has a winning $2$-tactic in
   $VSG([\kappa]^{<\lambda})$ for all $\kappa>\lambda$?
\end{problem}

   Our next theorem (Theorem \ref{3tacticth}) applies to abstract free ideals
   whose $\sigma$-completions have small principal bursting number.
   It is not clear to us whether ``3" occurring in Theorem
   \ref{3tacticth} is optimal. One of its applications is that ZFC+GCH
   implies that TWO has a winning 3-tactic in $VSG([\kappa]^{<\aleph_0})$
   for all $\kappa$. It is {\em very} likely that the ``3'' appearing
   in this application is not optimal, as will be discussed later.
 
\begin{theorem}\label{3tacticth}Let $J$ be a free ideal such that:
\begin{enumerate}
\item{$bu(\langle J\rangle,\subset) = \aleph_{2}$,}
\item{$add(\langle J\rangle,\subset) = \aleph_{1}$,}
\item{$cof(\langle J\rangle,\subset) = \lambda$ and}
\item{$[\lambda]^{<\aleph_0}$ has the coherent decomposition
property.}
\end{enumerate}
Then player TWO has a winning $3$-tactic in $VSG(J)$. 
\end{theorem}

\begin{description}\item[Proof.]{
   Let $J$ be a free ideal (on a set $S$) as in the hypotheses. Let
   ${\cal A}$ be a well-founded cofinal 
   family of cardinality $\lambda$, such that $|\{B\in{\cal
   A}:B\subseteq A\}|\leq\aleph_1$ for each $A\in{\cal A}$.  

   For each
   $A\in{\cal A}$ fix $\nu_A\leq\omega_{1}$ and a bijective
   enumeration $\{J_{\xi}(A):\xi < \nu_A\}$ of the set $\{X\in
   {\cal A}:X\subseteq A\}$.

   Choose a sequence $(C_{\xi}:\xi<\omega_{1})$ from $\langle J\rangle$ such
   that: 
\begin{enumerate}
\item{$C_{\xi}\subset C_{\nu}$ for $\xi < \nu$ and}
\item{$\cup_{\xi<\omega_{1}}C_{\xi}\not\in \langle J\rangle$.}
\end{enumerate} 

   For $A\in{\cal A}$ define
$\xi_A=\min\{\xi<\omega_1:C_{\xi}\not\subseteq A\}$.

   For $A\subset B$ elements from ${\cal A}$, define a set $\tau(A,B)$ such
   that $(S_{1},\dots,S_{n})$ is in $\tau(A,B)$ if:
\begin{enumerate}
\item{$2\leq n <\omega$,}
\item{$S_{1}=B$ and $S_{2}=A$,}
\item{$S_{j+1}\in\{J_{\xi}(S_{j}):\xi<\nu_{S_j}\mbox{ and }
C_{\xi}\subset S_{j-1}\}$ for $2\leq j< n$.} 
\end{enumerate}
   For $(S_{1},\dots,S_{n})$ and $(T_{1},\dots,T_{m})$ in $\tau(A,B)$ define
   $(S_{1},\dots,S_{n}) < (T_{1},\dots,T_{m})$ if $n<m$ and
   $(S_{1},\dots,S_{n})=(T_{1},\dots,T_{n})$. Then $(\tau(A,B),<)$ is a tree.
   Each branch of this tree is finite since $({\cal A},\subset)$ is
   well-founded. Indeed, $\tau(A,B)$ is a countable set.

   Define $F(A,B)$ to be the set of $X\in {\cal A}$ such that $X\in
   \{S_{1},\dots,S_{m}\}$ for some $(S_{1},\dots,S_{m})\in \tau(A,B)$. Then
   $F(A,B)$ is a countable set. Notice that if $C\subset A\subset B$
   are elements of ${\cal A}$ such that $C\in
   \{J_{\xi}(A):\xi\leq\nu_A \mbox{ and } C_{\xi}\subset B\}$,
   then $F(C,A)\subset F(A,B)$.  

   Let ${\cal B}\subset [{\cal A}]^{\aleph_0}$ be cofinal,
   well-founded and with the coherent decomposition property. For each
   $B\in{\cal B}$ choose a decomposition $B=\cup_{n=1}^{\infty}B^n$ where
   each $B^n$ is finite, and these decompositions satisfy the coherent
   decomposition requirement. By Proposition 15 of \cite{S2} we also fix
   a function 
$${\cal K}:[{\cal B}]^2\rightarrow\omega$$
   which witnesses that $({\cal
   B},\subset)\not\rightarrow(\omega-path)^2_{\omega/<\omega}$. 

   Define $\Phi_1:[{\cal A}]^2\rightarrow{\cal B}$ such that 

$$(\cup\{F(X,Y):(\exists(S_1,\dots,S_n)\in\tau(A,B))(X\subset Y \mbox{
   and }X,Y\in\{S_1,\dots,S_n\})\})
$$

   is a subset of $\Phi_1(A,B).$   Also define $\Phi_2:[{\cal
   A}]^2\rightarrow{\cal A}$ such that 

$$C_{\xi}\cup C_{\xi_B}\cup(\cup\Phi_1(A,B))\subset\Phi_2(A,B)$$
  
   where now $A=J_{\xi}(B)$.

   Note that if $A, \ B,$ and $C$ are elements of ${\cal A}$ such that
   $A\subset B\subset \Phi_2(A,B)\subset C$, then
   $\Phi_1(A,B)\subset\Phi_1(B,C)$. 

   Finally, choose for each $A\in {\cal A}$ a $\Phi_3(A)\in{\cal A}$ such
   that $A\cup C_{\xi_A}\subseteq \Phi_3(A)$.   

   Choose for each $A\in {\cal A}$ a sequence of sets
   $A^{0}\subseteq\dots,A^{n}\subseteq\dots$ such that each $A^{i}$ is
   in $J$ and $A=\cup_{n=0}^{\infty}A^{n}$. 

   We now define a $3$-tactic for TWO: First note that for the very
   strong game we may make the harmless assumption that player ONE's
   moves are all from the cofinal family ${\cal A}$. Let $A\subset
   B\subset C$ be sets from ${\cal A}$. Here
   are player TWO's responses ${\cal F}(A)$, ${\cal F}(A,B)$ and
   ${\cal F}(A,B,C)$:\\

Case 1: ${\cal F}(A)=(\emptyset,\Phi_3(A))$\\

Case 2: ${\cal F}(A,B)=(\emptyset,\Phi_2(A,B))$\\

Case 3: ${\cal F}(A,B,C)=(D,\Phi_2(B,C))$\\
   if $\Phi_2(A,B)\subseteq C$, where $D=C^m_1\cup\dots\cup C^m_r$ is given by:
   $m\geq {\cal K}(\{\Phi_1(A,B),\Phi_1(B,C)\})$ is minimal such that
   $(\Phi_1(A,B))^n\subseteq(\Phi_1(B,C))^n$ for all $n\geq m$, and
   $(\Phi_1(B,C))^m=\{C_1,\dots,C_r\}$. \\

Case 4: In all other cases define ${\cal F}(A,B,C)={\cal
   F}(B,C)$.\\

   To see that ${\cal F}$ is a winning $3$-tactic for TWO, consider a play
   $$(O_{1},(T_{1},S_{1}),O_{2},(T_{2},S_{2}),\dots)$$ of $VSG(J)$ for which
\begin{enumerate}
\item{$(T_{1},S_{1})={\cal F}(O_{1})$,}
\item{$(T_{2},S_{2})={\cal F}(O_{1},O_{2})$ and}
\item{$(T_{n+3},S_{n+3})={\cal F}(O_{n+1},O_{n+2},O_{n+3})$}
\end{enumerate}
   for all $n$.

   Then $T_1=T_2=\emptyset$, $S_1=\Phi_3(O_1)$, $S_2=\Phi_2(O_1,O_2)$
   and $S_{n+1}=\Phi_2(O_n,O_{n+1})$ for all $n\geq 2$. From the fact that
   $O_n\supseteq S_{n-1}$ for all $n\geq 2$ it follows that
$$ O_1\subset O_2\subset\Phi_2(O_1,O_2)\subseteq O_3\subset
   \Phi_2(O_2,O_3)\subseteq O_4\subset\dots,
$$
   whence
   $\Phi_1(O_1,O_2)\subset\Phi_1(O_2,O_3)\subset\Phi_1(O_3,O_4)\subset\dots$.
   For each $k$ let $m_k$ denote the minimal integer such that
\begin{enumerate}
\item{${\cal K}(\{\Phi_1(O_k,O_{k+1}),\Phi_1(O_{k+1},O_{k+2})\})\leq
   m_k$ and}
\item{$(\Phi_1(O_k,O_{k+1}))^n\subseteq(\Phi(O_{k+1},O_{k+2}))^n$ for
   all $n\geq m_k$.}
\end{enumerate}
   From the properties of ${\cal K}$ it follows that there are
   infinitely many $k$ such that $m_j<m_k$ for each $j<k$. Fix $i$, and
   fix the smallest $j\geq i$ such that $O_i\in\Phi_1(O_j,O_{j+1})$. Then
   let $t$ be minimal such that $O_i\in(\Phi_1(O_j,O_{j+1}))^t$. Then for
   each $k$ such that $m_{\ell}<m_k$ for each $\ell<k$, and $t<m_k$,
   $O^{m_k}_i\subseteq T_k$. It follows that
   $O_i\subseteq\cup_{n=1}^{\infty}T_n$. From this it follows that TWO
   won this ${\cal F}$-play of $VSG(J)$. $\QED$}
\end{description}

\begin{corollary}[GCH]\label{vsg:countablefinite}
For every infinite cardinal number $\kappa$, TWO has a winning
$3$-tactic in $VSG([\kappa]^{<\aleph_0})$.
\end{corollary}

   The results of Corollaries \ref{cor:countablecof4} and
\ref{vsg:countablefinite}  should be compared
   with those of Koszmider \cite{Ko} for the game
   $MG([\kappa]^{<\aleph_{0}})$. In Corollary \ref{cor:countablecof4}
   we show that there is a proper class of $\kappa$ such that TWO has a
   winning $2$-tactic in $SMG([\kappa]^{<\aleph_0})$, and thus in
   $VSG([\kappa]^{<\aleph_0})$. This class includes
   $\aleph_n$ for all $n<\omega$. In \cite{Ko} it is proven that TWO has a
   winning 2-tactic in $MG([\aleph_{n}]^{<\aleph_{0}})$ for 
   all $n\in \omega$ (\cite{Ko}, Theorem 18). Under the additional set
   theoretic assumption that both
   $\square_{\lambda}$ holds and $\lambda^{\aleph_{0}}=\lambda^{+}$ for all
   uncountable cardinal numbers $\lambda$ which are of countable cofinality,
   Koszmider further proves that player TWO has a winning 2-tactic in
   $MG([\kappa]^{<\aleph_{0}})$ for all $\kappa$ (\cite{Ko}, Theorem 19). In
   light of these results it is consistent that TWO 
   has a winning 2-tactic in the game $SMG([\kappa]^{<\aleph_{0}})$
   and thus in the game $VSG([\kappa]^{<\aleph_0})$ for all
   $\kappa$. 

    All this evidence leads us to believe that one could prove
   (without recourse to additional set theoretic hypotheses) that
   player TWO has a winning 2-tactic in the game
   $VSG([\kappa]^{<\aleph_{0}})$ for all infinite $\kappa$. We suspect
   even more: that TWO has a winning 2-tactic in
   $SMG([\kappa]^{<\aleph_0})$ for all $\kappa$. We state this
   formally as a conjecture:   
\begin{conjecture}\label{conj:countable-finite} Player TWO has a winning
   2-tactic in the game $SMG([\kappa]^{<\aleph_{0}})$ for each infinite
   cardinal number $\kappa$. 
\end{conjecture}

   One can modify the proof of Theorem \ref{3tacticth} to obtain the
   following result:

\begin{theorem}\label{3tacticthgeneralized}Let $J$ be a free ideal on a set
   $S$ such that 
\begin{enumerate}
\item{$bu(\langle J\rangle,\subset)=\aleph_n$ for some finite $n$,}
\item{there is an $(\omega_k,\omega_k)$-pseudo Lusin set in $(\langle
   J\rangle,\subset)$ for each $k\in\{1,\ldots,n\}$,}
\item{$cof(\langle J\rangle,\subset)=\lambda)$, and}
\item{$([\lambda]^{<\aleph_0},\subset)$ has the coherent decomposition
   property.}
\end{enumerate}
   Then player TWO has a winning $n+1$-tactic in $VSG(J)$.
\end{theorem}

   We now give an example which shows, assuming the Continuum
   Hypothesis, that the hypothesis that $add(\langle J\rangle,\subset)= 
   \aleph_1$ of Theorem \ref{3tacticth} is necessary (see Corollary
   \ref{cor:counterexample}). 

\begin{theorem}Let $\omega_{\alpha}$ be the initial ordinal
   corresponding with ${\frak c}$. Then there is a free ideal
   $J\subset{\EuScript P}(\omega_{\alpha+1})$ such that $cof(\langle
   J\rangle,\subset)=\aleph_{\alpha+1}$ and there is no positive integer $k$ for which
   TWO has a winning $k$-tactic in $VSG(J)$.
\end{theorem}
\begin{description}\item[Proof]{Define $J\subset{\EuScript
   P}(\omega_{\alpha+1})$ such that 
   $X\in J$ if, and only if, $|X|\leq\aleph_{\alpha}$ and
   $X\cap\omega$ is finite. 
   Then $cof(\langle J\rangle,\subset)=add(\langle J\rangle,\subset)=
   \omega_{\alpha+1}.$ By Theorem \ref{SMG,VSG} it suffices to show
   that TWO doesn't have a winning $2$-tactic in $SMG(J)$. \\

   Let $F$ be a $2$-tactic for TWO in $SMG(J)$. For
   $\omega<\eta<\omega_{\alpha+1}$ put 
   $\phi(\eta)=\sup(\eta\cup F(\eta)).$ Let
   $C\subseteq\omega_{\alpha+1}\backslash(\omega+1)$ be a closed
   unbounded set such that 
   $\phi(\gamma)<\beta$ whenever $\gamma<\beta$ are in $C$.\\

   For each $\eta\in C$ define
   $\phi_{\eta}:C\backslash(\eta+1)\rightarrow \omega_{\alpha+1}$ so that
   $\phi_{\eta}(\beta)=\sup(\beta\cup F(\eta,\beta))$ for all $\beta$. Then
   choose a closed, unbounded set $C_{\eta}\subseteq C\backslash(\alpha+1)$
   such that $\phi_{\eta}(\beta)<\gamma$ whenever $\beta<\gamma$ are in
   $C_{\eta}$.

   Let $D$ be the diagonal intersection of $(C_{\eta}:\eta\in C)$; i.e.,
   $D=\{\xi\in C:\xi\in\cap\{C_{\eta}:\eta<\xi \mbox{ and }\eta\in C\}$.
   Then $D$ is an unbounded subset of $\omega_{\alpha+1}$. Now observe that if
   $\eta_1<\eta_2<\eta_3$ are elements of $D$, then
\begin{enumerate}
\item{$\eta_2\in C_{\eta_1}$,}
\item{$\eta_3\in C_{\eta_1}\cap C_{\eta_2}$, and thus}
\item{$F(\eta_1)\subseteq\eta_2$ and $F(\eta_1,\eta_2)\subseteq
   \eta_3$.}
\end{enumerate}

   Define $\Phi:[D]^{2}\rightarrow\omega$ so that
$$\Phi(\eta,\beta)=\max(\omega\cap(F(\eta)\cup F(\eta,\beta))).$$
   By the Erd\"os-Rado theorem we obtain an $n<\omega$ and an uncountable
   $X\subset D$ such that $\Phi(\eta,\beta)=n$ for all $\eta<\beta\in X$. 
   Pick $\eta_1<\eta_2<\dots<\eta_m<\dots$ from $X$ and put
   $O_n=\eta_n$ for each $n$. Put $T_1=F(O_1)$ and $T_{n+1}=F(O_n,O_{n+1})$
   for each $n$.

   Then $(O_1,T_1,\dots,O_n,T_n,\dots)$ is an $F$-play of $SMG(J)$ which is
   lost by TWO. $\QED$}
\end{description}

\begin{corollary}\label{cor:counterexample} Assume the Continuum Hypothesis.
   Then there is a free ideal $J\subset {\EuScript P}(\omega_2)$ such
   that $cof(\langle J\rangle,\subset)=\aleph_2$, and there is no
   positive integer $k$ for which $TWO$ has a winning $k$-tactic in
   $VSG(J)$.
\end{corollary}

   We don't know if there is for each $m$ a free ideal $J_m$ such that TWO
does not have a winning $m$-tactic, but does have a winning $m+1$-tactic in
$VSG(J_m)$.

\begin{problem}Is there for each $m$ a free ideal $J_m$ such that TWO does not
   have a winning $m$-tactic, but does have a winning $m+1$-tactic in
   $VSG(J_m)$? 
\end{problem}

\subsection{The Banach-Mazur game and an example of Debs}

   The Banach-Mazur game is defined as follows for a topological space
   $(X,\tau)$. Players ONE and TWO alternately choose nonempty open subsets
   from $X$; in the n-th inning player ONE first chooses $O_{n}$ and TWO
   responds with $T_{n}$. An inning is played for each positive integer. The
   sets chosen by the players must satisfy the rule
$$
O_{n+1}\subseteq T_{n}\subseteq O_{n}
$$
   for all n. Player TWO wins the play 
$$
(O_{1},T_{1},\dots,O_{n},T_{n},\dots)
$$
   if the intersection of these sets is nonempty; otherwise player ONE wins.
   Following Galvin and Telgarsky \cite{G-T}, we denote this game by
   $BM(X,\tau)$.
   In the early 1980's Debs \cite{D} solved Problem 3 of \cite{F-K} by giving
   examples of topological spaces $(X,\tau)$ for which player TWO has a winning
   strategy in the game $BM(X,\tau)$, but no winning $1$-tactic. In all but one
   of Debs' examples it was known (in $ZFC$) that TWO has a winning
   $2$-tactic. We show here that also for the remaining example player
   TWO has a winning $2$-tactic (Corollary \ref{vstobm}). This was
   previously known under the assumption of some additional hypotheses.
  
   This result eliminates this example as a candidate for providing
   evidence (consistent, modulo $ZFC$) towards the following
   conjecture of Telgarsky: 

\begin{conjecture}[Telgarsky, \cite{T}, p. 236] For each positive integer $k$
   there is a topological space $(X_{k},\tau_{k})$ such that TWO does
   not have a winning $k$-tactic, but does have a winning $k+1$-tactic
   in the game $BM(X_{k},\tau_{k})$. 
\end{conjecture}

   The following unpublished result of Galvin is the only theorem
   known to us which gives general conditions under which TWO has a
   winning $2$-tactic if TWO has a winning strategy in the Banach-Mazur
   game:

\begin{theorem}[Galvin, unpublished] Let $(X,\tau)$ be a topological
   space for which TWO has a winning strategy in the Banach-Mazur game.
   If this space has a pseudo base ${\EuScript P}$ with the property that
\begin{itemize}
\item{$|\{V\in{\EuScript P}:B\subseteq V\}|<s(B)$ for each $B$ in
   ${\EuScript P}$,}
\end{itemize}
   then TWO has a winning $2$-tactic.
\end{theorem}

   Here the cardinal number $s(B)$ is defined to be the minimal
   $\kappa$ such that $B$ does not contain a collection of $\kappa$
   pairwise disjoint nonempty open subsets; it is said to be the {\em
   Souslin number} of $B$.

   This subsection is organised as follows. We first prove a theorem
   concerning $k$-tactics in the Banach-Mazur game which is analogous
   to Theorem 
   5 of \cite{G-T}. It provides an equivalent formulation of Telgarsky's
   conjecture which allows player TWO slightly more information: TWO may also
   remember the inning number. After this we give our result on Debs'
   example.
\subsubsection{Markov $k$-tactics.}

   Whereas a $k$-tactic for player TWO remembers at most the latest $k$ moves of
   the opponent, a strategy for TWO which remembers in addition to this
   information also the number of the inning in progress is called a {\em
   Markov $k$-tactic}. This choice of terminology is by analogy with the
   terminology ``tactic" (used by Choquet \cite{C}, p. 116, Definition 7.11
   for what we call a $1$-tactic) and ``Markov strategy" (used by Galvin and
   Telgarsky \cite{G-T}, p. 52 for what we call a Markov $1$-tactic). A
   k-tactic is the special case of a Markov k-tactic where the inning number is
   ignored by the player.

   Note that if $(X,\tau)$ has a dense set of isolated points then player TWO
   has a winning $1$-tactic in $BM(X,\tau)$. Thus we may assume that if at all
   possible, player ONE will avoid playing an open set which contains an isolated
   point. From the point of view of $k$-tactics for TWO we may therefore
   restrict our attention to topological spaces without isolated points. By the
   following proposition we may further restrict our attention to topological
   spaces in which each nonempty open set contains infinitely many
   pairwise disjoint open subsets.

\begin{prop}\label{filterdecomp}Let $(X,\tau)$ be a topological space with no
   infinite set of pairwise disjoint open subsets. Then there is a positive
   integer $n$ such that: 
$$
\tau\backslash\{\emptyset\}=\tau_{1}\cup\dots\cup\tau_{n}
$$
   where each $\tau_{i}$ has the finite intersection property.
\end{prop}
\begin{description}\item[Proof.]{

Claim 1: There is a positive integer n such that every collection of pairwise
disjoint nonempty open subsets is of cardinality $\leq n$. 

Proof of Claim 1: This is a well known fact: see e.g. \cite{C-N}, Lemma 2.10,
   p. 31. $\QED$

   Now let $n$ be the minimal positive integer satisfying Claim 1. Let
   ${\cal U}=\{U_{1},\dots,U_{n}\}$ be a collection of pairwise disjoint nonempty
   open subsets of the space. Then ${\cal U}$ is
   a maximal pairwise disjoint family.

   For $1\leq i\leq n$, let $\tau_{i}$ be a maximal family of nonempty open
   sets such that:
\begin{enumerate}
\item{$U_{i}\in\tau_{i}$,}
\item{any two elements of $\tau_{i}$ have nonempty intersection.}
\end{enumerate}

   Claim 2: $\tau\backslash\{\emptyset\}=\tau_{1}\cup\dots\cup\tau_{n}$

   Proof of Claim 2: Assume the contrary and let $Y$ be a nonempty open set
   which is in none of the $\tau_{i}$. Then we find for each $i$ an $X_{i}$ in
   $\tau_{i}$ which is disjoint from $Y$ (by maximality of each $\tau_{i}$). We
   may assume that $X_{i}\subseteq U_{i}$ for each $i$. But then
   $\{X_{1},\dots,X_{n},Y\}$ is a collection of $n+1$ pairwise disjoint
   nonempty open subsets of $(X,\tau)$, contradicting the choice of $n$.
   $\QED$ 

	Each $\tau_{i}$ has the finite intersection property. $\QED$} 
 \end{description}

\begin{prop}\label{finitaryth}Let $(X,\tau)$ be a topological space for which:
\begin{enumerate}
\item{Player TWO has a winning strategy in the game $BM(X,\tau)$ and}
\item{every collection of pairwise disjoint open subsets is finite.}
\end{enumerate}
   Then TWO has a winning $1$-tactic in $BM(X,\tau)$.
\end{prop}
\begin{description}\item[Proof.]{Write, by Proposition \ref{filterdecomp},
$$
\tau\backslash\{\emptyset\}=\tau_{1}\cup\dots\cup\tau_{n}
$$
   where each $\tau_{i}$ has the finite intersection property, and $n$ is
   minimal. Choose a pairwise
   disjoint collection $\{U_{1},\dots,U_{n}\}$ such that $U_{j}\in \tau_{j}$
   for each $j$.

 Claim 3: For each $j$, if $S_{1}\supseteq S_{2}\supseteq\dots$ is a
   denumerable chain from $\tau_{j}$, then $\cap_{n=1}^{\infty}S_{n}\neq
   \emptyset$.

Proof of Claim 3: Assume the contrary, and fix $j$ and a chain $S_{1}\supseteq
   S_{2}\supseteq\dots$ in $\tau_{j}$ such that $\cap_{n=1}^{\infty}S_{n}=
   \emptyset$. We may assume that $S_{n+1}\subset S_{n}\subset U_{j}$ for all
   $n$. 

   Let $F$ be a winning perfect information strategy for TWO in $BM(X,\tau)$.
   Consider the play $$(O_{1},T_{1},\dots,O_{m},T_{m},\dots)$$ which is defined
so that:
\begin{enumerate}
\item{$O_{1}=S_{1}$,}
\item{$T_{m}=F(O_{1},\dots,O_{m})$ for all $m$ and}
\item{$O_{m+1}=T_{m}\cap S_{m+1}$.}
\end{enumerate}
   Note that since each $S_{m}$ is a subset of $U_{j}$, each response by player
   TWO using $F$ is a member of $\tau_{j}$, whence each $O_{m}$ is a legal move
   by ONE. We now get the contradiction that TWO lost this play
   despite the fact that TWO was playing according to a winning strategy. This
   completes the proof of Claim 3. $\QED$
   
   We now define a winning $1$-tactic, $G$, for TWO. Let $U$ be a nonempty open
   subset of $X$. Choose the minimal $j$ such that $U_{j}\cap U\neq\emptyset$
   and put $G(U)=U_{j}\cap U$. Claim 3 implies that this is a winning
   $1$-tactic for TWO. $\QED$}
\end{description}

  \begin{theorem}\label{reductionth}Let $k$ be a positive integer. If
   player TWO 
   has a winning Markov $k$-tactic in the Banach Mazur game on some topological
   space, then TWO has a winning $k$-tactic in the Banach-Mazur game on that
   space. 
\end{theorem}
\begin{description}\item[Proof.]{Let $k$ be a positive integer and let
   $(X,\tau)$ be a topological space such that TWO has a winning Markov
   $k$-tactic in the game $BM(X,\tau)$. We may assume that
   this space has no isolated points. By Proposition \ref{finitaryth} we may also
   assume that every nonempty open subset of $X$ contains infinitely many
   pairwise disjoint open subsets (player ONE may safely avoid playing open
   subsets not having this property). By Theorem 5 of \cite{G-T} we may
   assume that $k>1$.

   Let $F$ be a winning Markov $k$-tactic for TWO. For each nonempty open set
   $U$, let $\{J_{m}(U):0<m<\omega\}$ bijectively enumerate a collection of
   infinitely many pairwise disjoint nonempty open subsets of $U$.

   Define a $k$-tactic $G$ for TWO as follows. Let $U_{1}\supseteq\dots\supseteq
   U_{j}$ be nonempty open sets, where $1\leq j\leq k$.

Case 1: $j=1$: Put $G(U_{1})=F(J_{2}(U_{1}),1)$.

Case 2: $j>1$ and $U_{i+1}\subseteq J_{l+i+1}(U_{i})$ for $1\leq i<j$, for some
        $l$. Put
        $G(U_{1},\dots,U_{j})=F(J_{l+2}(U_{1}),\dots,J_{l+j+1}(U_{j}),l+j)$.

Case 3: In all other cases define $G(U_{1},\dots,U_{j})=G(U_{j})$.

   To see that $G$ is a winning $k$-tactic for TWO, consider a play 
$$
(O_{1},T_{1},\dots,O_{m},T_{m},\dots)
$$
   such that 
\begin{itemize}
\item{$T_{j}=G(O_{1},\dots,O_{j})$ for $j\leq k$ and}
\item{$T_{n+k}=G(O_{n+1},\dots,O_{n+k})$ for all $n$.}
\end{itemize}

   From the definition of $G$ and the rules of the Banach-Mazur game it follows
  that $T_{1}$ is defined by Case 1 and $T_{m}$ for $m>1$ by Case 2. In
  particular, writing $S_{n}$ for $J_{n+1}(O_{n})$ we find that:

\begin{enumerate}
\item{$T_{j}=F(S_{1},\dots,S_{j},j)$ for $j\leq k$ and}
\item{$T_{n+k}=F(S_{n+1},\dots,S_{n+k},n+k)$}
\end{enumerate}
   for all n. Indeed,
$$
O_{1}\supseteq S_{1}\supseteq T_{1}\supseteq O_{2}\supseteq S_{2}\supseteq\dots
$$

   Since $F$ is a winning Markov $k$-tactic, it follows that
   $\cap_{n=1}^{\infty}O_{n}\neq\emptyset$. $\QED$}
\end{description}

\subsubsection{Debs' example}
   Let $\sigma$ be the topology of the real line whose elements are of the
   form $U\backslash M$ where $U$ is open and $M$ is meager in the usual
   topology. The symbol $BM({\reals},\sigma)$ denotes the Banach-Mazur game,
   played on the topological space $({\reals},\sigma)$. It is known that TWO
   has a winning strategy but does not have a winning $1$-tactic in $BM({\reals
   },\sigma)$. 

\begin{corollary}\label{vstobm} Player TWO has a winning 2-tactic in
   the game $BM({\reals},\sigma)$.  
\end{corollary}

\begin{description}\item[Proof.]{Theorem 22 of \cite{S1} and Corollary
\ref{corsmg}. $\QED$} 
\end{description}

\section{Appendix: Consistency of the hypotheses of
   Theorem \protect\ref{noktactic}.}

   We start with a ground model $V$ and let ${\bold
   P} \in V$ be  a forcing notion of cardinality $\leq {\frak c}$. For a
   cardinal $\kappa$, denote by $\Pk$ the product of $\kappa$ copies
   of $\bold P$ taken side-by-side with countable supports.

\begin{lemma}\label{lemma:consistency2} Let $\lambda$ be an uncountable
   cardinal. Suppose:
\begin{enumerate}
\item{ $\kappa\geq\kappa_1\geq\kappa_2\geq\kappa_3\geq\omega_2$ are
cardinal numbers such that
\begin{itemize}
\item{$\kappa$ is a regular cardinal,}
\item{$\kappa\rightarrow(\kappa_1)^2_{\lambda}$,}
\item{$\kappa_1\rightarrow(\kappa_2)^3_{{\frak c}}$,}
\item{$\kappa_2\rightarrow(\kappa_3)^2_{\lambda}$ and}
\item{$\kappa_3\rightarrow(\omega_2)^3_{{\frak c}}$.}
\end{itemize}}
\item{Forcing with ${\bold P}$ adds a real to the ground model.} 
\end{enumerate}

   Then $\twoal \ropl$ holds in the forcing extension $\VPk$.
\end{lemma}

\begin{description}\item[Proof.]{ Let $\lambda, \kappa, \kappa_1,
   \dots , \kappa_3, \bold P$ be as in the 
   assumptions. 
   Our argument closely follows section 2 of \cite{To1}.

   For sets $A,B$ the symbol $A/B$ denotes $\{ \ab : \ \alpha \in A , 
   \ \beta \in B , \ \alpha < \beta \}$.

   Note that $\VPk$ satisfies $\twoal = \kappa$; we prove
   that $\kappa \ropl$ holds in $\VPk$.

   Let $[\kappa ]^2 = \bigcup_{i < \lambda} {\dot K}_i$ be a given partition
   in $\VPk$. Let $\dot U$ be a $\Pk$-name for a member of 
   $[\kappa ]^\kappa$. Pick $A \in [\kappa ]^\kappa$ and for each 
   $\alpha \in A$, a $q_\alpha \in \Pk$ such that 
   $q_\alpha \forces \alpha \in {\dot U}$ and such that the
   $q_\alpha$'s form a $\Delta$-system. Define $H: [A]^2 \rightarrow
   (\lambda+1)$ so that $H(\{ \alpha, \beta \}) = i$ if $i$ is the
   minimal $j$ such that $p \forces \{ \alpha , \beta \} \in {\dot
   K}_j$ for some $p \leq q_\alpha, 
   q_\beta$ if such $j$ exists (i.e., if $q_\alpha$ and $q_\beta$ are
   compatible), and $H(\{ \alpha , \beta \} ) = \lambda$ if $q_\alpha$
   is incompatible with $q_\beta$.

   By our choice of $\kappa$, the partition relation 
   $\kappa \rightarrow (\kappa_1 )^2_\lambda$ holds.
   Therefore, choose $A_1 \subset [A]^{\kappa_1}$
   and $i \leq \lambda$ such that 
   $H''[A_1]^2 = \{ i \}$. Since $\Pk$ satisfies the $\twoal^+$-c.c., we
   have $i < \lambda$.

   Let $\langle \pab : \ \{ \alpha , \beta \} \in [A_1]^2 \rangle$
   be a fixed sequence of conditions such that $\pab \leq q_\alpha , q_\beta$
   and $\pab \forces \ab \in \dKi$.
   For $\alpha < \beta < \gamma$ in $A_1$ we define 
   $H_0 (\{ \alpha , \beta , \gamma \} )$ to be a pair 
   $(c,d)$, where $c$ codes $\pab$ and $\pag$ as structures as well as
   relations  
   between the ordinals of $dom ( \pab)$ and $dom (\pag )$, and $d$
   does the same  
   for $\pag$ and $\pbg$. Since there are only $\twoal$ such pairs,
   and since $\kappa_1 \rightarrow (\kappa_2 )^3_{\twoal}$ holds, choose 
   $A_2 \in [A_1]^{\kappa_2}$ and $(c,d)$ such that 
   $H_0''[A_2]^3 = \{ (c,d) \}$. 
   For convenience, assume that $A_2$ has order type $\kappa_2$.
   It follows that for each 
   $\alpha \in A_2$ the sequence $\langle \pab : \ \beta \in A_2 \backslash
   (\alpha + 1) \rangle$ forms a $\Delta$-system with root $p^0_\alpha$
   ($\leq q_\alpha$), and that for each $\gamma \in A_2$ the sequence
   $\langle \pbg : \ \beta \in A_2 \cap \gamma \rangle$ forms a 
   $\Delta$-system with root $p^1_\gamma$ ($\leq q_\gamma$). Moreover, the
   $p_\alpha^0$'s and $p^1_\gamma$'s form $\Delta$-systems with roots
   $p^0$ and $p^1$ respectively. To see the latter, note that we may 
   shrink $A_2$ to a cofinal subset $A_3$ so that the relevant
   $p^0_\alpha$'s and $p^1_\alpha$'s do in fact form a $\Delta$-system.
   Now consider $\alpha , \beta , \gamma \in A_3$, and $\alpha ' , \beta '
   \in A_2$. Comparing $H_0(\{ \alpha , \beta , \gamma \})$, 
   $H_0( \{ \alpha , \beta ' , \gamma \} )$ and $H_0 ( \{\alpha ' ,
   \beta ' , \gamma \} )$, one sees that the sequence
   $\langle p_\alpha^0 : \ \alpha \in A_2 \rangle$ forms a
   $\Delta$-system. A similar argument works for the $p^1_\gamma$'s.

   Also, $p^0$ is compatible with $p^1$.
   We call $\langle \pab : \ \ab \in B/B \rangle$ a {\it double 
   $\Delta$-system} with root $p^0 \cup p^1$. 

   There is no reason why for a given $\alpha$ the conditions
   $p^0_\alpha$ and $p^1_\alpha$ should be compatible: if these were
   always compatible, our argument would yield a consistency proof of
   $\twoal \rightarrow (\omega )^2_\lambda$, which is false in $ZFC$.

   We now save as much of the compatibility between $p_\alpha^0$ and
   $p_\alpha^1$ as is needed for the consistency proof of
   $\twoal \ropl$. Thin out $A_2$ to a cofinal
   subset $A_3$ such that 
   $dom(p^0_\alpha \cup p^1_\alpha ) 
   \cap dom(p^0_\beta \cup p^1_\beta) = dom (p^0 \cup p^1)$ for all
   $\ab \in A_3 / A_3$. Then in particular 
   $p^1_\alpha$ and $p^0_\beta$ are compatible for $\ab \in A_3/A_3$.

   Now repeat the reasoning above with $A_2$ in place of $A$,
   $\kappa_2$ in place of $\kappa$, $\kappa_3$ in place of 
   $\kappa_1$, and $\omega_2$ in place of $\kappa_2$. Also, 
   $p^1_\alpha$ will now play the role of $q_\alpha$, and 
   $p^0_\beta $ the role of $q_\beta$ for $\{ \alpha , \beta \} \in 
   A_3 / A_3$. 
   We get a set $A_4 \subset A_3$ of order type $\omega_2$ and 
   some $j < \lambda$ (which may be different from $i$), conditions 
   $\bpab$ for $\ab \in A_4 / A_4$ that form a double $\Delta$-system
   with root $\bpz \cup \bpo$, and we get roots $\bpaz$ and $\bpgo$ as
   before. 
   Now  $\bpab \forces \ab \in {\dot K}_j$ for 
   $\ab \in A_4 / A_4$.

   Our choice of $\bpab$ at the beginning of the second run of the argument
   insures that $\bpaz \leq p_\alpha^1$ and $\bpgo \leq p^0_\gamma$, 
   and hence $\bpz \leq p^1$ and $\bpo \leq p^0$.

\smallskip

Now let $\bold G$ be a generic subset of $\Pk$. Define:

$\dX = \{ \alpha \in A_4 : \ p^0_\alpha \in {\bold G}\}$,

$\dY = \{ \alpha \in A_4 : \ p^1_\alpha \in {\bold G} \}$,

$\dW = \{ \alpha \in A_4 : \ \bpaz \in {\bold G} \}$,

$\dZ = \{ \alpha \in A_4 : \ \bp^1_\alpha \in {\bold G} \}$.

\smallskip

   Then $\dZ \subset \dX$ and $\dW \subset \dY$, and all four sets are
   cofinal in $A_4$.

   Now $\bp^0 \cup \bp^1$ forces the following facts:

(1) $\exists \delta_1 \in \omega_2 \forall \alpha \in \dX \backslash \delta_1
\ \{ \beta \in \dW : \ \ab \in \dKi \} $ is cofinal in $A_4$, and

(2) $\exists \delta_2 \in \omega_2 \forall \alpha \in \dY \backslash 
\delta_2 \ \{ \beta \in \dZ : \ \ab \in {\dot K}_j \}$ is cofinal in $A_4$.

\smallskip

   The combination of (1) and (2) suffices to construct in $\VPk$ an 
   $\omega$-path of the given partition that uses only colors $i$ and $j$:

   Let $\delta = max \{ \delta_1 , \delta_2 \}$. Inductively define
   an increasing sequence $\langle x_n: \ n \in \omega \rangle$ of
   ordinals such that $x_{2k} \in Z$ (and hence in $X$), 
   $x_{2k+1} \in W$, and $\{ x_{2k}, x_{2k+1} \} \in \dKi$ 
   (by (1)); $\{ x_{2k+1}, x_{2k+2} \} \in {\dot K}_j$ (by (2)).

   It remains to prove (1) and (2). We shall prove (1) only; the proof of (2)
   is similar, and is a special case of
   \cite{To1}, section 2, property (1)].

   Assume that $\bp^0 \cup \bp^1$ does not force (1). Then we can find 
   a condition $\bp^2 \leq \bp^0 \cup \bp^1$ and a $\Pk$-name 
   $\dot D \in [\dX ]^{\omega_2}$ and for each $\beta \in \dot D$ a
   $\gamma_\beta \in A_4 \backslash (\beta + 1)$ 
   such that $\bp^2  \forces  \{ \beta_ , \delta \}
   \notin \dKi$ whenever $\delta \in \dot W \backslash \gamma_\beta$.

   Working in $V$, we pick $B \in [A_4]^{\omega_2}$ such that for each 
   $\beta \in B$ we find $r_\beta \leq p^0_\beta \cup \bp^2$ 
   such that $r_\beta \forces \beta \in \dot D$, and $r_\beta$ decides
   the value of $\gamma_\beta$. We may assume that the $r_\beta$'s form
   a $\Delta$-system with root $\leq \bp^2 \leq \bp^0 \cup \bp^1$, and that 
   $\gamma_\beta < \delta$ for all $\{\beta , \delta\} \in B/B$.
   Since $\langle \pab : \ \beta \in B \backslash (\alpha + 1) \rangle$ forms a
   $\Delta$-system, we may also assume that 
   $dom(r_\beta) \cap dom(p_{\beta  , \delta}\backslash p^0_{\beta} ) 
   = \emptyset $ for all $ \delta > \gamma_\beta$ in $A_4$.

   Pick $\delta \in A_2$ such that $B \cap \delta$ is uncountable and
   $dom(\bp^0_\delta ) \cap dom (\bp^2) = dom (\bp^0)$. Since 
   $\langle p_{\beta , \delta} : \ \beta \in B \cap \delta \rangle$ forms a 
   $\Delta$-system with root $p^1_\delta$ and since $dom(\bp^0_\delta )$ 
   is countable, we have
   $dom(p_{\beta , \delta} \backslash p^1_\delta) \cap dom (\bp^0_\delta ) 
   \not= \emptyset$ for only countably many $\beta \in B \cap \delta$.
   So pick a $\beta \in B \cap \delta$ such that 
   $dom(p_{\beta , \delta} \backslash p^1_\delta) \cap dom (\bp^0_\delta ) 
   = \emptyset$.
\smallskip

Define $r \in \Pk$ as follows:

\smallskip

$dom (r) = dom (r_\beta ) \cup dom (\bp^0_\delta ) \cup
dom (p_{\beta , \delta } \backslash p^1_\delta )$,

\smallskip

$r| dom (r_\beta \cup \bp^0_\delta ) = r_\beta \cup \bp^0_\delta$,

\smallskip

and 

\smallskip

$r(\xi ) = p_{\beta , \delta } (\xi ) $ for $\xi \in 
dom (p_{\beta , \delta } \backslash dom (r_\beta \cup \bp^0_\delta ))$.

\smallskip

   Then $r$ is a well-defined condition with the property that
   $r \leq r_\beta , \bp^0_\delta$ and $p_{\beta , \delta}$. 
   So $r$ forces that $\{\beta , \delta\} \in X/W$  and 
   that $\{ \beta , \delta \} \in K_i$, which is a contradiction. $\QED$}
\end{description}

   If $\Pk$ is as 
   in the assumptions of Lemma \ref{lemma:consistency2}, then $\Pk$ is
   a $\twoal^+$-c.c. poset.
   If GCH holds in the ground model and $\lambda = \omega_1$, then 
   our proof works if $\kappa \geq \aleph_8$. One can obtain the
   consistency of $\twoal \rightarrow (\omega -path)^2_{\omega_1/<\omega}$  
   with a smaller size of the continuum, but this is not essential for
   our purposes.  Todorcevic has for example shown that, adjoining at
   least $\omega_2$ 
   Cohen reals to a model of the Continuum Hypothesis, produces a model
   in which $\omega_2\rightarrow(\mbox{$\omega$-path})^2_{\omega/<3}$.\\

   We have actually proved something apparently stronger than 
   $\twoal \ropl$ in $\VPk$, namely a relation denoted by
   $\twoal \rightarrow (\omega - path)^2_{\lambda/<3}$.

   We do not know an answer to the following two
   problems concerning the $\omega$-path partition relation:

\begin{problem} Is it for each integer $k>2$ consistent, for some
   infinite cardinal numbers $\kappa$ and $\lambda$, that
   $\kappa\not\rightarrow(\omega\mbox{-path})^{2}_{\lambda/<k}$, but  
   $\kappa\rightarrow(\omega\mbox{-path})^{2}_{\lambda/<k+1}$?
\end{problem}

\begin{problem} Is it consistent, for some infinite
   cardinal numbers $\kappa$ and $\lambda$, that for each $k<\omega$,
   $\kappa\not\rightarrow(\omega\mbox{-path})^{2}_{\lambda/<k}$, but  
   $\kappa\rightarrow(\omega\mbox{-path})^{2}_{\lambda/<\omega}$?
\end{problem}

\begin{theorem}[Todorcevic]\label{consistency2} If $ZFC$ is a
   consistent theory, then so is the theory $ZFC \ + \ cof(\langle\JR\rangle ,
   \subset ) = \aleph_1 \ + \  {\frak c} \ropo$.
\end{theorem}

\begin{description}\item[Proof]{
   Theorem \ref{consistency2} is an immediate consequence of Lemma
   \ref{lemma:consistency2}: It is well known that  if CH holds in the ground
   model, and $\bold P$ is e.g. Sacks or Prikry-Silver forcing, then (b)
   and (c) of 
   the lemma hold for every $\kappa$. It is also known that adding any
   number of Sacks or Prikry-Silver reals side-by-side with countable
   supports to a model of CH, one obtains a model where the collection
   of meager sets whose Borel codes are from the ground model, is a
   cofinal subfamily of $\JR$ (see \cite{M}). Since $| ^{\omega}\omega\cap
   V| = \aleph_1$, we get $cof(\langle\JR\rangle , \subset ) =
   \aleph_1$ in the forcing extension. $\QED$} 
\end{description}

\end{document}